\input ./style/arxiv-general.cfg
\documentclass[aos,MSNbibl,nameyear,seceqn,rotating,dvips]{arximspdf}
\makeatletter
   \@ifpackageloaded{graphicx}{}{\usepackage{graphicx}}
\makeatother
\usepackage{accents}
\usepackage[mathscr]{eucal}
\usepackage{mathbh,multirow}

\doi{10.1214/15-AOS1355}
\volume{43}
\issue{6}
\pubyear{2015}
\firstpage{2653}
\lastpage{2675}
\docsubty{FLA}

\makeatletter
\newcommand{\rrvert}{\vert}
\newcommand{\rrVert}{\Vert}
\newcommand{\llvert}{\vert}
\newcommand{\llVert}{\Vert}
\newproclaim{newremark}{Remark}[section]
\newtheorem{theorem}{Theorem}[section]

\newtheorem{corollary}[theorem]{Corollary}
\newproclaim{definition}[theorem]{Definition}
\newproclaim{exmp}[example]{Example}
\newproclaim{rmrk}[remark]{Remark}

\def\eqdef{\stackrel{\operatorname{def}}{=}}

\def\tow{\stackrel{w}{\longrightarrow}}

\newcommand{\cc}[1]{\mathscr{#1}}
\newcommand{\bb}[1]{\mathbf{#1}}
\newcommand{\bbb}[1]{\bolds{#1}}

\renewcommand{\bar}[1]{\overline{#1}}
\renewcommand{\hat}[1]{\widehat{#1}}
\renewcommand{\tilde}[1]{\widetilde{#1}}

\def\Var{\operatorname{Var}}

\def\argmax{\mathop{\operatorname{argmax}}}

\def\tr{\operatorname{tr}}

\def\R{\mathbb{R}}
\def\E{\mathbb{E}}
\def\P{\mathbb{P}}

\def\T{\top}

\def\dv{\bb{d}}

\def\Yv{\bb{Y}}

\def\gammav{\bbb{\gamma}}

\def\xiv{\bbb{\xi}}

\def\CONST{\mathtt{C}}

\def\rddelta{\delta}

\def\Ind{\mathbh{1}}

\def\DP{D}
\def\DPc{\DP_{0}}

\def\gmi{\mathtt{b}}

\def\xivb{\xiv_{\rdb}}

\def\ex{\mathrm{e}}

\def\gm{\mathtt{g}}

\def\yy{\mathtt{y}}

\def\xx{\mathtt{x}}

\def\gmu{\mathfrak{a}}

\def\rups{\rr_{0}}

\def\CS{\cc{E}}

\def\nunu{\nu_{0}}

\def\LL{\cc{L}}

\def\dimp{p}

\def\thetav{\bbb{\theta}}
\def\thetavs{\thetav^{*}}

\def\thetavb{\thetav^{\dag}}

\def\Thetas{\Theta_{0}}

\def\VP{V}
\def\VPc{\VP_{0}}

\def\rr{\mathtt{r}}

\def\zz{\mathfrak{z}}
\def\zzb{\tilde{\zz}}
\def\tt{\mathfrak{t}}

\def\smb{\operatorname{smb}}

\def\dv{h}

\def\DP{D}
\def\DPc{\DP_{0}}
\def\sbt{\circ}

\def\sbtt{\circ}

\def\Lb{L^{\sbt}}
\def\Eb{\E^{\sbt}}
\def\Pb{\P^{\sbt}}
\def\Varb{\Var^{\sbt}}

\def\xivb{\xiv^{\sbt}}

\def\zzb{\zz^{\sbt}}

\def\deltasmb{\delta_{\smb}}

\def\xivbr{\bar{\xiv}}

\def\xivbbr{{\xivbr}^{\sbt}}

\def\Bsmb{B_{0}}

\def\DeltaBEfull{{\Delta}_{{\operatorname{B.E.}}, \full}}
\def\DeltaBEfullzz{{\Delta}_{{\operatorname{B.E.}}\ \zz, \full}}

\def\DeltabfullI{{\Delta}_{\mb, \full}}

\def\CONST{\mathtt{C}}

\def\DeltaW{\Delta_{\operatorname{W}}}

\def\ZZ{\mbox{$\mathfrak{Z}$}}

\def\Deltaqqf{\Delta_{\mathtt{qf},1}}
\def\Deltaqlf{\Delta_{\mathtt{qf},2}}

\def\full{\operatorname{full}}

\def\mb{\operatorname{b}}
\def\sm{\operatorname{sm}}

\def\Id{\mathbf{I}}

\def\Hc{H_{0}}

\def\thetavt{\tilde{\thetav}}
\def\thetavb{\tilde{\thetav}}
\def\thetavbt{\tilde{\thetav}^{\sbtt}}

\def\dpc{d_{0}}

\newcommand{\eqref}[1]{(\ref{#1})}

\def\NDi{\cc{N}}
\def\yy{\xx}
\def\ED{\mathit{E D}}
\def\SD{\mathit{S D}}
\def\SmB{\operatorname{SmB}}
\makeatother

\begin{document}
\begin{frontmatter}

\title{Bootstrap confidence sets under model misspecification}
\runtitle{Bootstrap confidence sets}

\begin{aug}
\author[A]{\fnms{Vladimir} \snm{Spokoiny}\thanksref{M1,M2,M3,M4,M5,T1}\ead[label=e1]{spokoiny@wias-berlin.de}}
\and
\author[A]{\fnms{Mayya} \snm{Zhilova}\corref{}\thanksref{M1,T2}\ead[label=e2]{zhilova@wias-berlin.de}}
\runauthor{V. Spokoiny and M. Zhilova}
\thankstext{T1}{Supported in part by
Laboratory for Structural Methods of Data Analysis in Predictive
Modeling, MIPT,
RF government Grant, ag. 11.G34.31.0073,
by the German Research Foundation (DFG) through the
Research Unit 1735 and
supported in part by the Russian Science Foundation grant (project 14-50-00150).}
\thankstext{T2}{Supported by the German Research Foundation (DFG)
through the Collaborative
Research Center 649 ``Economic Risk'' and
supported in part by the Russian Science Foundation grant (project
14-50-00150).}

\affiliation{Weierstrass Institute for Applied Analysis and
Stochastics\thanksmark{M1},
Humboldt University Berlin\thanksmark{M2},
Moscow Institute of
Physics and Technology\thanksmark{M3},
Institute for Information Transmission Problems RAS\thanksmark{M4} and
Higher School of Economics, Moscow\thanksmark{M5}}

\address[A]{Weierstrass-Institute\\
Mohrenstr. 39\\
10117 Berlin\\
Germany\\
\printead{e1}\\
\phantom{E-mail:\ }\printead*{e2}}

\end{aug}

%
\received{\smonth{11} \syear{2014}}
%
\revised{\smonth{6} \syear{2015}}


\begin{abstract}
A multiplier bootstrap procedure for construction of likelihood-based
confidence sets is considered for finite samples and a possible model
misspecification.
Theoretical results justify the bootstrap validity for a small or
moderate sample size
and allow to control the impact of the parameter dimension $ p $:
the bootstrap approximation works if $ p^{3}/n $ is small.
The main result about bootstrap validity
continues to apply even if the underlying parametric model is misspecified
under the so-called small modelling bias condition.
In the case when the true model deviates significantly from the
considered parametric family, the bootstrap procedure is still
applicable but it becomes a bit conservative:
the size of the constructed confidence sets is increased by the
modelling bias.
We illustrate the results with numerical examples for misspecified
linear and logistic regressions.
\end{abstract}
%
%

\end{frontmatter}

\section{Introduction}
Since introducing in 1979 by \citet{Efron1979}, the bootstrap procedure
became one of the most powerful
and common tools in statistical confidence estimation and hypothesis testing.
Many versions and extensions of the original bootstrap method have been
proposed in the literature;
see, for example, \citet{Barbe1995weighted},
\citet{Bucher2013multiplier},
\citet{ChattBose2005generalized},
\citeauthor{Chen2009efficient}
(\citeyear{Chen2009efficient,Chen2014sieve}),
\citet{Horowitz2001bootstrapHandbook},
\citet{Janssen1994weighted},
\citet{Lavergne2013smooth},
\citet{MaKosorok2005robust},
\citet{Mammen1993bootstrap},
\citet{NewtonRaft1994WLB},
\citet{Wu1986wildboot}
among many others.
This paper focuses on the multiplier bootstrap procedure which
attracted a lot of attention last time
due to its nice theoretical properties and numerical performance.
We mention the papers of \citet{ChattBose2005generalized}, \citet
{ArlotBlanch2010} and \citet{ChernoMultBoot} for the most advanced
recent results.
\citet{ChattBose2005generalized} showed some results on asymptotic
bootstrap consistency in a very general framework for estimators
obtained by solving estimating equations.
\citet{ChernoMultBoot} presented a number of nonasymptotic results on
bootstrap validity with applications to special
problems like testing many moment restrictions or parameter choice for
a LASSO procedure.
\citet{ArlotBlanch2010} constructed a nonasymptotical confidence bound
in $\ell_{s}$-norm ($s\in[1,\infty]$) for the mean of a sample of
high dimensional i.i.d. Gaussian vectors (or with a symmetric and
bounded distribution), using the generalized weighted bootstrap for
resampling of the quantiles.

This paper makes a further step in studying the multiplier bootstrap
method in the problem of
confidence estimation by a quasi maximum likelihood method.
For a rather general parametric model, we consider likelihood-based
confidence sets
with the radius determined by a multiplier bootstrap.
The aim of the study is to check the validity of the bootstrap
procedure in situations with a growing parameter dimension, a limited
sample size, and
a possible misspecification of the parametric assumption.
The main result of the paper explicitly describes the error term of the
bootstrap approximation.
This particularly allows to track the impact of the parameter dimension
$ \dimp$
and of the sample size $ n $ in the quality of the bootstrap procedure.
As one of the corollaries, we show bootstrap validity under the constraint
``$ \dimp^{3}/n $-small.''
\citet{ChattBose2005generalized} stated results under the condition ``$
\dimp/n $-small''
but their results only apply to low dimensional projections of the MLE vector.
In the likelihood-based approach, the construction involves the
Euclidean norm of the MLE
which leads to completely different tools and results.
\citet{ChernoMultBoot} allowed a huge parameter dimension with ``$ \log
(\dimp) / n $ small''
but they essentially work with
a family of univariate tests which again differs essentially from the
maximum likelihood
approach.

Another interesting and important issue is the impact of the model
misspecification
on the accuracy of bootstrap approximation.
A surprising corollary of our error bounds is that
the bootstrap confidence set
can be used even if the underlying parametric model is slightly misspecified
under the so-called \emph{small modelling bias} $(\SmB)$ condition.
If the modelling bias becomes large, the bootstrap confidence sets are
still applicable,
but they become more and more conservative.
$(\SmB)$ condition is given in Section~\ref{sect_conditions} and
it is consistent with classical
bias--variance relation in nonparametric estimation.

Our theoretical study uses the square-root Wilks (sq-Wilks) expansion
from \citeauthor{Sp2012Pa}
(\citeyear{Sp2012Pa,Spokoiny2013Bernstein}) which approximates
the square root likelihood ratio statistic
by the norm of the standardized score vector.
Further, we extend the sq-Wilks expansion to the bootstrap log-likelihood
and adopt the Gaussian approximation theory (GAR) to the special case
when the distribution
of the Euclidean norm of a non-Gaussian vector is approximated by the
distribution
of the norm of a Gaussian one with the same first and second moments.
The Gaussian comparison technique based on the Pinsker inequality completes
the study and allows to bridge the real unknown coverage probability
and the conditional bootstrap
coverage probability under $(\SmB)$ condition.
In the case of a large modelling bias, we state a one-sided bound: the
bootstrap quantiles are uniformly larger than the real ones.
This effect is nicely confirmed by our simulation study.

Now consider the problem and the approach in more detail.
Let the data sample $\Yv= (Y_{1},\ldots, Y_{n} )^{\T}$
consist of \textit{independent} random observations and belong to the
probability space $ (\Omega, \mathcal{F},\P )$.
We do not assume that the observations $Y_{i}$ are identically
distributed; moreover, no specific parametric structure of $\P$ is
being required.
In order to explain the idea of the approach we start here with a
parametric case, however, assumption \eqref{assum:PA} below is not
required for the results. Let $\P$ belong to some known regular
parametric family $ \{\P_{\thetav} \}\eqdef \{ \P_{\thetav
}\ll\mu_{0}, \thetav\in\Theta\subset\R^{\dimp} \}$. In this
case, the true parameter $\thetavs\in\Theta$ is such that
%
\begin{equation}
\label{assum:PA} \P\equiv\P_{\thetavs}\in\{\P_{\thetav}\},
\end{equation}
and the initial problem of finding the properties of unknown
distribution $\P$ is reduced to the equivalent problem for the
finite-dimensional parameter $\thetavs$.
The parametric family $\{\P_{\thetav}\}$ induces the log-likelihood
process $L(\thetav)$ of the sample $\Yv$,
\[
L(\thetav)=L(\Yv,\thetav)\eqdef\log{ \biggl(\frac{d \P_{\thetav}}{d
\mu_{0}}(\Yv) \biggr)}
\]
and the maximum likelihood estimate (MLE) of $\thetavs$,
%
\begin{equation}
\label{def:thetavt} \tilde{\thetav}\eqdef\argmax _{\thetav\in\Theta} L(\thetav).
\end{equation}
The asymptotic Wilks phenomenon [\citet{wilks1938}] states
that for the case of i.i.d. observations with the sample size tending
to the infinity the likelihood ratio statistic converges in
distribution to $\chi^{2}_{\dimp}/2$, where $\dimp$ is the
parameter dimension
\[
2 \bigl\{L(\thetavt)-L\bigl(\thetavs\bigr) \bigr\}\tow\chi^{2}_{\dimp},\qquad
n\to \infty.
\]
Define the likelihood-based confidence set as
%
\begin{equation}
\label{CS} \CS(\zz)\eqdef \bigl\{\thetav: L(\thetavt)-L(\thetav)\leq
\zz^{2}/2 \bigr\},
\end{equation}
then the Wilks phenomenon implies
\[
\P \bigl\{\thetavs\in\CS(\zz_{\alpha, \chi^{2}_{p}}) \bigr\}\to\alpha ,\qquad n\to\infty,
\]
where $\zz_{\alpha, \chi^{2}_{p}}^{2}$ is the $(1-\alpha)$-quantile for the $\chi^{2}_{\dimp}$ distribution.
This result is very important and useful under the parametric
assumption, that is, when \eqref{assum:PA} holds.
In this case, the limit distribution of the likelihood ratio is
independent of the model parameters or in other words it is \textit
{pivotal}. By this result, a sufficiently large sample size allows to
construct the confidence sets for $\thetavs$ with a given coverage
probability.
However, a possibly low speed of convergence of the likelihood ratio
statistic makes the asymptotic Wilks result hardly applicable to the
case of small or moderate samples.
Moreover, the asymptotical pivotality breaks down if the parametric
assumption \eqref{assum:PA} does not hold [see \citet{huber1967}] and,
therefore, the whole approach may be misleading if the model is
considerably misspecified.
If the assumption \eqref{assum:PA} does not hold, then the ``true''
parameter is defined by the projection of the true measure $\P$ on
the parametric family $\{\P_{\thetav}\}$:
%
\begin{equation}
\label{thetavs} \thetavs \eqdef \argmax _{\thetav\in\Theta} \E L(\thetav).
\end{equation}
The recent results by \citeauthor{Sp2012Pa}
(\citeyear{Sp2012Pa,Spokoiny2013Bernstein})
provide a nonasymptotic version of square-root Wilks phenomenon for
the case of misspecified model.
It holds with an exponentially high probability
%
\begin{equation}
\label{Wilks_intro} \bigl\llvert \sqrt{2 \bigl\{L(\thetavt)-L\bigl(\thetavs\bigr) \bigr\}}-\|
\xiv\| \bigr\rrvert \leq \DeltaW\simeq\frac{\dimp}{\sqrt{n}},
\end{equation}
where $\xiv\eqdef\DPc^{-1}\nabla_{\thetav}L(\thetavs)$,
$\DPc^{2}\eqdef-\nabla_{\thetav}^{2}\E L(\thetavs)$.
The bound is nonasymptotical, the approximation error term $\DeltaW$
has an explicit form (the precise statement is given in Theorem~B.2,
Section~B.1 of the supplementary material [\citet
{SpZh2014PMBsupp}], and it depends on the parameter dimension $\dimp$,
sample size $n$ and the probability of the random set on which the
result holds.

Due to this bound,
the original problem of finding a quantile of the LR test statistic $
L(\thetavt)-L(\thetavs) $
is reduced to a similar question for the approximating quantity $\|\xiv
\|$.
The difficulty here is that in general $\|\xiv\|$ is nonpivotal,
it depends on the unknown distribution $\P$ and the target parameter
$\thetavs$.

In the present work, we study the \textit{multiplier bootstrap} (or
\textit{weighted bootstrap}) procedure for estimation of the quantiles
of the likelihood ratio\vadjust{\goodbreak} statistic. The idea of the procedure is to
mimic a distribution of the likelihood ratio statistic by reweighing
its summands with random multipliers independent of the data
\[
\Lb(\thetav) \eqdef \sum_{i=1}^{n}\log{
\biggl(\frac{d\P_{\thetav}}{d\mu
_{0}}(Y_{i}) \biggr)}u_{i}.
\]
Here, the probability distribution is taken conditionally on the data $
\Yv$, which is denoted by the sign $^{\sbt}$ (also $\Eb$ and $
\Varb$ denote expectation and variance operators w.r.t. the
probability measure conditional on $\Yv$). The random weights $
u_{1},\ldots,u_{n}$ are i.i.d., independent of $\Yv$ and it holds for them:
$\Eb(u_{i})=1$, $\Varb(u_{i})=1$, $\Eb\exp(u_{i})<\infty$.
Therefore, the multiplier bootstrap induces the probability space
conditional on the data $\Yv$.
A simple but important observation is that
$ \Eb\Lb(\thetav) \equiv L(\thetav) $, and hence,
\[
\argmax _{\thetav} \Eb\Lb(\thetav) = \argmax _{\thetav} L(\thetav) =
\tilde{\thetav}. \label{bEbLttt}
\]
This means that the target parameter in the bootstrap world is
precisely known and it coincides with
the maximum likelihood estimator
$\thetavt$ conditioned on $\Yv$, therefore, the bootstrap
likelihood ratio statistic  $\Lb(\thetavbt)-\Lb(\thetavt)\eqdef\sup_{\thetav\in\Theta}\Lb(\thetav)-\Lb(\thetavt)$ is fully computable
and leads to a simple computational procedure for the approximation of
the distribution of $L(\thetavt)-L(\thetavs)$.

The goal of the present study is to show in a nonasymptotic way the
validity of the described multiplier bootstrap procedure and to obtain
an explicit bound on the error of coverage probability.
In other words, we are interested in nonasymptotic approximation of
the distribution of
$ \{L(\thetavt)-L(\thetavs) \}^{1/2}$\vspace*{1pt} with the distribution of
$ \{\Lb(\thetavbt)-\Lb(\thetavb) \}^{1/2}$. So far there
exist very few theoretical nonasymptotic results about bootstrap validity.
Classical asymptotic tools for showing the bootstrap consistency are based
on weak convergence arguments which are not applicable in the finite
sample set-up.
Some different methods have to be applied.
In particular, the approach of \citet{Liu1988} based on Berry--Esseen theorem
can be extended to a finite sample set-up with a univariate parameter.
For a high dimensional parameter space, important contributions are
done in the
recent papers by \citet{ArlotBlanch2010} and \citet{ChernoMultBoot}.
The latter paper used a Gaussian approximation, Gaussian comparison and
Gaussian anti-concentration technique in high dimension.
Our approach is similar but we combine it with the square-root Wilks expansion
and use Pinsker's inequality for Gaussian comparison and
anti-concentration steps.
The main steps of our theoretical study are illustrated by the
following scheme:\looseness=1
%
%
\begin{eqnarray}\label{rectangle}\qquad
\tabcolsep=0pt
\begin{tabular}{lllllll}
&& \multicolumn{1}{c}{sq-Wilks} &&\multicolumn{1}{c}{Gauss.}&&\\
&& \multicolumn{1}{c}{theorem}&&\multicolumn{1}{c}{approx.}&&\\
$\Yv$-world:& ${\sqrt{ 2L(\thetavt)-2L\bigl(\thetavs\bigr)}}$&
$\displaystyle\mathop{\approx}_{\dimp/\sqrt{n}}^{\mathrm{}}$& $\|\xiv\|$&
$\displaystyle\mathop{\mathop{\approx}_{ (\dimp^{3}/n
)^{1/8}}^{w}}$& $\|\xivbr\|$ &\\
&&&&&
$w \ {\rotatebox{-90}{$\approx$}}\ {{\sqrt{\dimp}\deltasmb^{2}}}$&
\multicolumn{1}{l}{\multirow{2}{37pt}[10pt]{\centering Gauss. compar.}}\\
\\
\multicolumn{1}{l}{\multirow{2}{43pt}[10pt]{\centering Bootstrap world:}}& ${\sqrt{2\Lb\bigl(\thetavbt\bigr)-2\Lb(\thetavt) }}$&
$\displaystyle\mathop{\approx}_{\dimp/\sqrt{n}}$& $\bigl\|\xivb\bigr\|$&
$ \displaystyle\mathop{\approx}^{w}_ {(\dimp^{3}/n
)^{1/8}} $& $\bigl\|{\xivbr}^{\sbt}\bigr\|$
\end{tabular}
\end{eqnarray}\looseness=0
where
\[
\xivb \eqdef \xivb\bigl(\thetavs\bigr) \eqdef\DPc^{-1}\nabla_{\thetav}
\bigl[\Lb\bigl(\thetavs\bigr)-\Eb \Lb\bigl(\thetavs\bigr) \bigr]. \label{xibootdef}
\]
The vectors $ \xivbr$ and $ \xivbbr$ are zero mean Gaussian and
they mimic the covariance structure
of the vectors $ \xiv$ and $ \xivb$:
$\xivbr\sim\NDi(0,\Var\xiv)$, $\xivbbr\sim\NDi (0,\Varb\xivb
)$.

The error term shown below each arrow corresponds to the i.i.d. case considered
in details in Section~\ref{typical_local}.
The upper line of the scheme corresponds to the $\Yv$-world, the
lower line---to the bootstrap world.
In both lines, we apply two steps for approximating the corresponding
likelihood ratio statistics.
The first approximating step is the nonasymptotic square-root Wilks
theorem: the bound \eqref{Wilks_intro} for the $\Yv$-case and a
similar statement for the bootstrap world, which is obtained in
Theorem~B.4, Section~B.2 in \citet
{SpZh2014PMBsupp}.
The corresponding error is of order $\dimp/\sqrt{n}$ for the case of
i.i.d. observations; in the bootstrap world the square-root Wilks
expansion implies
\[
\bigl\llvert \sqrt{2\Lb\bigl(\thetavbt\bigr)-2\Lb(\thetavt) } - \bigl\|{\xivb}(\thetavt)\bigr\|
\bigr\rrvert \leq\CONST\dimp/\sqrt{n} \label{sqWilksboot}
\]
for $\xivb(\thetav) \eqdef\DPc^{-1}\nabla_{\thetav} [\Lb(\thetav
)-\Eb\Lb(\thetav) ]$.
In our approximation diagram, we use $\xivb(\thetavs)$ instead of
$\xivb(\thetavt)$ which is more convenient for the GAR step and is
justified by
Lemma~B.7 in \citet{SpZh2014PMBsupp} showing that
$ \|\xivb(\thetavt)-\xivb(\thetavs) \|\leq\CONST\dimp/\sqrt
{n}$.

The next step is called \emph{Gaussian approximation} (GAR)
which means that the distribution of the Euclidean norm $\|\xiv\|$ of
a centered random vector
$ \xiv$ is close to the distribution of the similar norm of a
Gaussian vector $\|\xivbr\|$
with the same covariance matrix as $\xiv$.
A similar statement holds for the vector $\xivb$.
Thus, the initial problem of comparing the distributions of the
likelihood ratio statistics is reduced to the comparison of the
distributions of the Euclidean norms of two centered normal vectors $
\xivbr$ and $\xivbbr$ (Gaussian comparison).
This last step links their distributions and encloses the approximating scheme.
The Gaussian comparison step is done by computing the Kullback--Leibler
divergence between two multivariate Gaussian distributions [i.e., by
comparison of the covariance matrices of $\nabla_{\thetav}L(\thetavs
)$ and $\nabla_{\thetav}\Lb(\thetavs)$] and applying Pinsker's
inequality [Lemma~A.7 in \citet{SpZh2014PMBsupp}].
At this point, we need to introduce the ``small modelling bias''
condition $(\SmB)$ from Section~\ref{sect:ConditAddBoot}.
It is formulated in terms of the following nonnegative-definite $\dimp
\times\dimp$ symmetric matrices:
%
\begin{eqnarray}
\label{def:Hc} \Hc^{2}&\eqdef& \sum_{i=1}^{n}
\E \bigl[\nabla_{\thetav}\ell_{i}\bigl(\thetavs\bigr)\nabla _{\thetav}
\ell_{i}\bigl(\thetavs\bigr)^{\T} \bigr],
\\
\label{def:Bsmb} \Bsmb^{2}&\eqdef& \sum_{i=1}^{n}
\E \bigl[\nabla_{\thetav}\ell_{i}\bigl(\thetavs\bigr) \bigr]\E \bigl[
\nabla_{\thetav
}\ell_{i}\bigl(\thetavs\bigr) \bigr]^{\T}
\end{eqnarray}
for $\ell_{i}(\thetav)\eqdef\log (\frac{d \P_{\thetav}}{d \mu
_{0}}(Y_{i}) )$,
so that $\Var \{\nabla_{\thetav}L(\thetavs) \}=\Hc^{2}-\Bsmb
^{2}$.
If the parametric assumption \eqref{assum:PA} is true or if the data $
\Yv$ are i.i.d., then it holds $\E [\nabla_{\thetav}\ell
_{i}(\thetavs) ]\equiv0$ and $\Bsmb^{2}=0$.
The $(\SmB)$ condition roughly means that the bias term $ \Bsmb
^{2} $ is small relative to $ \Hc^{2} $.
Below we show that the Kullback--Leibler distance between
the distributions of two Gaussian vectors $ \xivbr$ and $ \xivbbr$
is bounded by $\dimp\|\Hc^{-1}\Bsmb^{2}\Hc^{-1}\|^{2}/2$.
The $(\SmB)$ condition precisely means that this quantity is
small [in scheme \eqref{rectangle} it is denoted by $\sqrt{\dimp
}\deltasmb^{2}$].
In Section~\ref{sect:smb_condit}, the value $\|\Hc^{-1}\Bsmb^{2}\Hc
^{-1}\|$ is evaluated for some commonly used models: the case of
i.i.d. observations, generalized linear model and linear quantile regression.
Below we distinguish between two situations: when the condition $(\SmB)$ is
fulfilled and the opposite case. Theorems~\ref
{thm:cumulat} and \ref{thm:cumulat2} in Section~\ref{sect_mainres} deal
with the first case. It provides the cumulative error term for the
coverage probability of the confidence set \eqref{CS}, taken at the $
(1-\alpha)$-quantile computed with the multiplier bootstrap procedure.
The proof of this result [see Section~B.4 in \citet
{SpZh2014PMBsupp}] summarizes the steps of scheme~\eqref{rectangle}.
The biggest term in the full error is induced by Gaussian approximation
and requires the ratio $\dimp^{3}/n$ to be small. In the case of a
``large modelling bias,'' that is, when $(\SmB)$ does not hold,
the multiplier bootstrap procedure continues to apply. It turns out
that the bootstrap quantiles increase with the growing modelling bias;
hence, the confidence set based on it remains valid, however, it may
become conservative. This result is given in Theorem~\ref{thm:resbias}
of Section~\ref{sect_mainres}.
The problems of Gaussian approximation and comparison for the Euclidean
norm are considered in Sections~A.2 and A.4 of the supplementary material [\citet{SpZh2014PMBsupp}]
in general terms independently of the statistical setting of the paper,
and might be interesting by themselves. Section~A.4
in \citet{SpZh2014PMBsupp} presents also an anti-concentration
inequality for the Euclidean norm of a Gaussian vector.
This inequality shows how the deviation probability changes with a threshold.
The general results on GAR are summarized in Theorem~A.1 in the
supplementary material [\citet
{SpZh2014PMBsupp}] and restated in Proposition~B.12 in \citet{SpZh2014PMBsupp} for the setting of
scheme \eqref{rectangle}.
These results are also nonasymptotic with explicit errors and apply
under the condition that
the ratio $\dimp^{3}/n$ to be small.

In Theorem~\ref{thm:res_BE}, we consider the case of a scalar parameter
$\dimp=1$ with an improved error term. Furthermore, in Section~\ref{sect:res_smooth} we propose a modified version of a quantile function
based on a smoothed probability distribution. In this case, the
obtained error term is also better than in the general result.

Notation: $\|\cdot\|$ denotes Euclidean norm for vectors and spectral
norm for matrices; $\CONST$ is a generic constant. The value $\xx
>0$ describes our tolerance level:
all the results will be valid on a random set of probability ($1- C \ex
^{-\yy}$)
for an explicit constant $C$.
Everywhere we give explicit error bounds and show how they depend on $
\dimp$ and $n$ for the case of the i.i.d. observations $Y_{1},\ldots
,Y_{n}$ and $\yy\leq\CONST\log{n}$. More details on it are given
in Section~\ref{typical_local}. In Section~B.3 in
the supplementary material [\citet{SpZh2014PMBsupp}], we also consider generalized
linear model and linear quantile regression, and show for them the
dependence on $\dimp$ and $n$ of all the values appearing in main
results and their conditions.

The paper is organized as follows: the main results are stated in
Section~\ref{sect_mainres}. Their proofs are given in
Sections~B.4, B.5 and B.6 of the supplementary material
[\citet{SpZh2014PMBsupp}].
Section~\ref{sect:numer} contains numerical results for misspecified
linear and logistic regressions. In Section~\ref{sect_conditions}, we
give all the required conditions, provide information about dependence
of the involved terms on $n$ and $\dimp$ and consider the $(\SmB)$ condition for some models.
Section~A in \citet{SpZh2014PMBsupp} collects some
useful statements on Gaussian approximation and
Gaussian comparison.

\section{Multiplier bootstrap procedure}
\label{sect_mainres}
Let $\ell_{i}(\thetav)$ denote the parametric log-density of
the $i$th observation
\[
\ell_{i}(\thetav)\eqdef\log \biggl(\frac{d \P_{\thetav}}{d \mu
_{0}}(Y_{i})
\biggr),
\]
then $L(\thetav)=\sum_{i=1}^{n} \ell_{i}(\thetav)$.
Consider i.i.d. scalar random variables $u_{i}$ independent of $\Yv
$ with $\E u_{i}=1$, $\Var u_{i}=1$, $\E\exp(u_{i})<\infty$ for
all $i=1,\ldots,n$.
Multiply the summands of the likelihood function $L(\thetav)$ with
the new random variables
\[
\Lb(\thetav)\eqdef\sum_{i=1}^{n}
\ell_{i}(\thetav)u_{i},
\]
then it holds
$\Eb\Lb(\thetav)=L(\thetav)$,
where $ \Eb$ stands for the conditional expectation given $ \Yv$.
Therefore, the quasi MLE for the $\Yv$-world is a target parameter
for the bootstrap world:
\[
\argmax _{\thetav\in\Theta}\Eb\Lb(\thetav) = \argmax _{\thetav\in\Theta}L(\thetav) =
\tilde{\thetav}.
\]
The corresponding quasi MLE under the conditional measure $\Pb$ is
defined as
\[
\thetavbt\eqdef\argmax _{\thetav\in\Theta}\Lb(\thetav).
\]
The likelihood ratio statistic in the bootstrap world is equal to $\Lb
(\thetavbt)-\Lb(\thetavt)$ in which all the entries are known including
the function $\Lb(\thetav)$ and the arguments $\thetavbt$, $
\thetavt$.

Let $1-\alpha\in(0,1)$ be an unknown desirable confidence level of
the set $\CS(\zz)$:
%
\begin{equation}
\label{CSineq} \P \bigl(\thetavs\in\CS(\zz) \bigr) \geq 1-\alpha.
\end{equation}
Here, the parameter $\zz\geq0$ determines the size of the confidence set.
Define $\zz_{\alpha}$ as
the minimal possible value of $ \zz$ such that \eqref{CSineq} is fulfilled:
%
\begin{equation}
\label{def:zzalpha} \zz_{\alpha} \eqdef \inf \bigl\{\zz\geq0 \colon \P
\bigl(L(\thetavt) - L\bigl(\thetavs\bigr) > \zz^{2}/2 \bigr) \leq\alpha \bigr\}.
\end{equation}
For evaluating this value, we apply the multiplier bootstrap procedure
which replaces the unknown data distribution with the artificial bootstrap
distribution given the observed sample.
The target value $ \zz_{\alpha} $ is approximated by
the value $\zzb_{\alpha}$ defined as the upper $\alpha$-quantile of
$  \{2\Lb(\thetavbt) - 2\Lb(\thetavt) \}^{1/2} $:
%
\begin{equation}
\label{def:zzbalpha} \zzb_{\alpha} \eqdef \inf \bigl\{\zz\geq0 \colon \Pb
\bigl(\Lb\bigl(\thetavbt\bigr) - \Lb(\thetavt) > \zz^{2}/2 \bigr) \leq \alpha
\bigr\}.
\end{equation}
Note that the bootstrap probability $ \Pb$ and log-likelihood excess
$ \Lb(\thetavbt) - \Lb(\thetavt) $
depends on the data $ \Yv$ and thus, $ \zzb_{\alpha} $ is random
as well.
Theoretical results of the next section justify the proposed approach.

\subsection{Main results}
Now we state the main results for the general set-up. The approximating
error terms and the conditions are specified in Section~B.3 of
the supplementary material [\citet{SpZh2014PMBsupp}] for
popular examples including i.i.d. observations, generalized regression
model and linear quantile regression.
Our first result claims that the random quantity
$ \Pb (\Lb(\thetavbt) - \Lb(\thetavt) > \zz^{2}/2  ) $ is
close in
probability to the value $ \P (L(\thetavt) - L(\thetavs) > \zz
^{2}/2  ) $
for a wide range of $ \zz$-values.

\begin{theorem}
\label{thm:cumulat}
Let the conditions of Section~\ref{sect_conditions} be fulfilled, then
it holds
for $\zz\geq\max\{2,\sqrt{\dimp}\}+\CONST(\dimp+\yy)/\sqrt{n}$
with probability $\geq1-12\ex^{-\yy}$:
\[
\bigl\llvert \P \bigl( L(\thetavt) - L\bigl(\thetavs\bigr) > \zz^{2}/2 \bigr) -
\Pb \bigl( \Lb\bigl(\thetavbt\bigr) - \Lb(\thetavt) > \zz^{2}/2 \bigr) \bigr
\rrvert  \leq \Delta_{\full}.
\]
The error term $\Delta_{\full} \leq\CONST\{(\dimp+\yy)^{3}/n\}
^{1/8}$ in the case
of i.i.d. model; see Section~\ref{typical_local}. Explicit definition
of the error term $\Delta_{\full}$ is given in Section~\textup{B.4}
of the supplementary material [\citet{SpZh2014PMBsupp}]; see \textup{(B.41)}
and \textup{(B.42)} therein.
\end{theorem}

The term $\Delta_{\full}$ can be viewed as a sum of the error terms
corresponding to each step in the scheme \eqref{rectangle}.
The largest error term equal to $\CONST\{(\dimp+\yy)^{3}/n\}^{1/8}$
is induced by GAR. This error rate is not always optimal for GAR, for
example, in the case of $\dimp=1$ or for the i.i.d. observations [see
Remark~A.2 in \citet{SpZh2014PMBsupp}]. In Theorems \ref
{thm:res_BE} and \ref{thm:smooth}, the rate is $\CONST\{(\dimp+\yy
)^{3}/n\}^{1/2}$.

The next result can be viewed as ``bootstrap validity.''

\begin{theorem} [(Validity of the bootstrap under a small modelling bias)]
\label{thm:cumulat2}
Assume the conditions of Theorem~\ref{thm:cumulat}.
Then for $\alpha\leq1-8\ex^{-\yy}$, it holds
\[
\bigl\llvert \P \bigl( L(\thetavt) - L\bigl(\thetavs\bigr) > \bigl(\zzb_{\alpha}\bigr)^{2}/2
\bigr) - \alpha \bigr\rrvert  \leq \Delta_{\zz, \full}. \label{PLtttsDzzf}
\]
The error term $ \Delta_{\zz, \full} \leq\CONST\{(\dimp+\yy)^{3}/n\}
^{1/8}$ in the case
of the i.i.d. model; see Section~\ref{typical_local}.
For a precise description, see \textup{(B.46)} and \textup{(B.47)} of the
supplementary material [\citet{SpZh2014PMBsupp}].
\end{theorem}

In view of definition \eqref{CS} of the likelihood-based confidence
set, Theorem~\ref{thm:cumulat} implies the following:

\begin{corollary}[(Coverage probability error)]
Under the conditions of Theorem~\ref{thm:cumulat2}, it holds that
\[
\bigl\llvert \P \bigl\{\thetavs\in\CS (\zzb_{\alpha} ) \bigr\} -(1-
\alpha)\bigr\rrvert \leq\Delta_{\zz, \full}.
\]
\end{corollary}

\begin{newremark}[(Critical dimension)]
\label{rem:applicab}
The error term $\Delta_{\full}$ depends on the ratio $ \dimp^{3}/n $.
The bootstrap validity can be only stated if this ratio is small.
The obtained error bound seems to be mainly of theoretical interest,
because the condition ``$ (\dimp^{3}/n)^{1/8} $ is small'' may
require a huge sample.
However, it provides some qualitative information about the bootstrap
behavior as the parameter
dimension grows.
Our numerical results show that the accuracy of bootstrap approximation
is very reasonable
in a variety of examples with $ \dimp\ll n $.
\end{newremark}

In the following theorem, we consider the case of the scalar parameter
\mbox{$\dimp=1$}. The obtained error rate is $1/\sqrt{n}$, which is
sharper than $1/n^{1/8}$. Instead of the GAR for the Euclidean norm
from Section~A in \citet{SpZh2014PMBsupp}, we use here
the Berry--Esseen theorem [see also Remark~A.2 in \citet
{SpZh2014PMBsupp}].

\begin{theorem}[(The case of $\dimp=1$, using the Berry--Esseen theorem)]
\label{thm:res_BE}
Let the conditions of Section~\ref{sect_conditions} be fulfilled.
\begin{longlist}[1.]
\item[1.]%
For $\zz\geq1+\CONST(1+\yy)/\sqrt{n}$, it holds with probability $
\geq1-12\ex^{-\xx} $
\[
\bigl\llvert \P \bigl( L(\thetavt) - L\bigl(\thetavs\bigr) > \zz^{2}/2 \bigr)-
\Pb \bigl( \Lb\bigl(\thetavbt\bigr) - \Lb(\thetavt) > \zz^{2}/2 \bigr)\bigr
\rrvert \leq \DeltaBEfull.
\]

\item[2.]%
For $ \alpha\leq1-8\ex^{-\yy}$
\[
\bigl\llvert \P \bigl( L(\thetavt) - L\bigl(\thetavs\bigr) > \bigl(\zzb_{\alpha}\bigr)^{2}/2
\bigr)-\alpha\bigr\rrvert \leq \DeltaBEfullzz.
\]
\end{longlist}
The error terms $\DeltaBEfull, \DeltaBEfullzz\leq\CONST(1+\yy)/\sqrt
{n}$ in the case \ref{typical_local}.
Explicit definitions of $\DeltaBEfull$ is given in \textup{(B.48) }and
\textup{(B.49)} in Section~\textup{B.4} of the supplementary material [\citet{SpZh2014PMBsupp}].
\end{theorem}

\begin{newremark}[(Bootstrap validity and weak convergence)]
The standard way of proving the bootstrap validity is based on weak
convergence arguments;
see, for example, \citet{Mammen1992does},
\citet{VaartWellner1996weak},
\citet{Janssen2003bootstrap}, \citet{ChattBose2005generalized}.
If the statistic $ L(\thetavt) - L(\thetavs) $ weakly converges to a
$ \chi^{2} $-type distribution,
one can state an asymptotic version of the results of Theorems \ref
{thm:cumulat}, \ref{thm:res_BE}.
Our way is based on a kind of nonasymptotic Gaussian approximation and
Gaussian comparison
for random vectors and allows to get explicit error terms.
\end{newremark}

\begin{newremark}[(Use of Edgeworth expansion)]
The classical results on confidence sets for the mean of population
states the accuracy of order
$ 1/n $ based on the second-order Edgeworth expansion; see \citet
{Hall1992bootstbook}. Unfortunately, if the considered parametric model
can be misspecified, even the leading term
is affected by the modelling bias, and the use of Edgeworth expansion
cannot help
in improving the bootstrap accuracy.
\end{newremark}

\begin{newremark}[(Choice of the weights)]
In our construction, similarly to \citet{ChattBose2005generalized},
we apply a general distribution of the bootstrap weights
$ u_{i} $ under some moment conditions.
One particularly can use Gaussian multipliers as suggested by \citet
{ChernoMultBoot}.
This leads to the exact Gaussian distribution of the vectors $ \xivb
$ and
is helpful to avoid one step of Gaussian approximation for these vectors.
\end{newremark}

\begin{newremark}[(Skipping the Gaussian approximation step)]
The biggest error term $\CONST\{(\dimp+\yy)^{3}/n\}^{1/8}$ in Theorem~\ref{thm:cumulat}
is induced by the Gaussian approximation step.
In some particular cases, the Gaussian approximation step can be
avoided leading to better error bounds.
For example, if the marginal score vectors $\nabla_{\thetav}\ell
_{i}(\thetavs)$ are normally distributed, and the random bootstrap
weights are normal as well,
$u_{i} \sim\NDi(1,1)$, then the vectors $\xiv$ and $\xivb$ are
automatically normal, and the GAR step can be skipped.
If the marginal score vectors $\nabla_{\thetav}\ell_{i}(\thetavs)$
are i.i.d. and symmetrically distributed [s.t. $\nabla_{\thetav}\ell
_{i}(\thetavs)\sim-\nabla_{\thetav}\ell_{i}(\thetavs)$], and the
centered bootstrap weights follow the Rademacher distribution [$u_{i}
\sim2\operatorname{Bernoulli}(0.5)$], then the recent results by \citet
{ArlotBlanch2010} can be applied to
show that the conditional distribution of
$\|\xivb(\thetavs)\|$ given the data is close to the distribution of
$\|\xiv\|$.
However, such methods require some special structural conditions on the
underlying measure $ \P$
like symmetricity or Gaussianity of the errors and may fail if these conditions
are violated.
It remains a challenging question how a nice performance of
a general bootstrap procedure even for small or moderate samples can be
explained.
\end{newremark}

Now we discuss the impact of modelling bias, which comes from a
possible misspecification
of the parametric model.
As explained by the approximating diagram \eqref{rectangle}, the
distance between
the distributions of the likelihood ratio statistics can be
characterized via the distance between two multivariate normal distributions.
To state the result, let us recall the definition of the full Fisher
information matrix $\DPc^{2}\eqdef-\nabla_{\thetav}^{2}\E L(\thetavs)$. For the matrices $\Hc^{2}$ and $\Bsmb^{2}$, given in \eqref
{def:Hc} and \eqref{def:Bsmb}, it holds $\Hc^{2}>\Bsmb^{2}\geq0$. If
the parametric assumption \eqref{assum:PA} is true or in the case of an
i.i.d. sample $\Yv$, $\Bsmb^{2}=0$. Under the condition $(\SmB)$ $\|\Hc^{-1} \Bsmb^{2}\Hc^{-1}\|$ enters linearly in the
error term $\Delta_{\full}$ in
Theorem~\ref{thm:cumulat}.

The first statement in Theorem~\ref{thm:resbias} below says that the
effective coverage probability of the confidence set based on the
multiplier bootstrap is \textit{larger} than the nominal coverage
probability up to the error term $\DeltabfullI\leq\CONST\{(\dimp+\yy
)^{3}/n\}^{1/8}$. The inequalities in the second part of Theorem~\ref
{thm:resbias} prove the \textit{conservativeness of the bootstrap
quantiles}: the quantity $\sqrt{\tr\{\DPc^{-1}\Hc^{2}\DPc^{-1}\}}-\sqrt
{\tr\{\DPc^{-1}(\Hc^{2}-\Bsmb^{2})\DPc^{-1}\}}\geq0$ increases with
the growing modelling bias.

\begin{theorem}
[(Performance of the bootstrap for a large modelling bias)]
\label{thm:resbias}
Under the conditions of Section~\ref{sect_conditions} except for $(\SmB)$, it holds for $\zz\geq\break \max\{2,\sqrt{\dimp}\}+\CONST(\dimp
+\yy)/\sqrt{n}$ with probability $\geq1-14\ex^{-\xx}$:
\begin{longlist}
\item[1.]
\[
\P\bigl (L(\thetavt)-L\bigl(\thetavs\bigr) >\zz^{2}/2 \bigr)
\leq
\Pb \bigl(\Lb\bigl(\thetavbt\bigr)-\Lb(\thetavt) >\zz^{2}/2 \bigr)+\DeltabfullI.
\]

\item[2.]
\begin{eqnarray*}
\zzb_{\alpha}&\geq&
\zz_{(\alpha+\DeltabfullI)}
\\
&&{}+ \sqrt{\tr\bigl\{\DPc^{-1}\Hc^{2}\DPc^{-1}\bigr\}}-\sqrt{\tr\bigl\{\DPc^{-1}\bigl(\Hc
^{2}-\Bsmb^{2}\bigr)\DPc^{-1}\bigr\}}
- \Deltaqqf,\\
\zzb_{\alpha}&\leq&
\zz_{(\alpha-\DeltabfullI)}
\\
&&{}+ \sqrt{\tr\bigl\{\DPc^{-1}\Hc^{2}\DPc^{-1}\bigr\}}-\sqrt{\tr\bigl\{\DPc^{-1}\bigl(\Hc
^{2}-\Bsmb^{2}\bigr)\DPc^{-1}\bigr\}}
+ \Deltaqlf.
\end{eqnarray*}
The term $\DeltabfullI\leq\CONST\{(\dimp+\yy)^{3}/n\}^{1/8}$ is
given in \textup{(B.51)} in Section~\textup{B.5} of the
supplementary material [\citet{SpZh2014PMBsupp}].
The positive values $\Deltaqqf, \Deltaqlf$ are given in \textup{(B.55)},
\textup{(B.54)} in Section~\textup{B.5} in \citet{SpZh2014PMBsupp}; they are bounded from
above with $(\gmu^{2}+\gmu_{B}^{2})(\sqrt{8\yy\dimp}+6\yy)$ for the
constants $\gmu^{2}>0,  \gmu_{B}^{2}\geq0$ from conditions $(\mathcal{I})$, $(\mathcal{I}_{B})$.
\end{longlist}
\end{theorem}

\begin{newremark}
There exists some literature on robust (and heteroscedasticity robust)
bootstrap procedures; see, for example,
\citet{Mammen1993bootstrap}, \citet{Aerts2001missp}, \citet
{Kline2012higher}. However, to our knowledge there are no robust
bootstrap procedures for the likelihood ratio statistic,
most of the results compare the distribution of the estimator obtained
from estimating equations,
or Wald/score test statistics with their bootstrap counterparts in the
i.i.d. setup.
In our context, this would correspond to the noise misspecification in
the log-likelihood function and it is addressed automatically by the
multiplier bootstrap.
Our notion of modelling bias includes the situation when the target
value $ \thetavs$
from \eqref{thetavs} only defines a projection (the best parametric
fit) of the data distribution.
In particularly, the quantities $ \E\nabla_{\thetav}\ell_{i}(\thetavs
) $ for different $ i $
do not necessarily vanish yielding a significant modelling bias.
Similar notion of misspecification is used in the literature on
Generalized Method of Moments;
see, for example, \citet{Hall2005GMMbook}.
Chapter~5 therein considers the hypothesis testing problem with two
kinds of misspecification: local and nonlocal, which would correspond
to our small and large modelling bias cases.

An interesting message of Theorem~\ref{thm:resbias} is that the
multiplier bootstrap
procedure ensures a prescribed coverage level for this target value $
\thetavs$ even without small modelling bias restriction;
however, in this case, the method is somehow conservative because the
modelling bias is transferred into
the additional variance in the bootstrap world. The numerical
experiments in Section~\ref{sect:numer} agree with this result.
\end{newremark}

\subsection{Smoothed version of a quantile function}
\label{sect:res_smooth}

This section explains how to improve the accuracy of bootstrap
approximation using a smoothed quantile function.
The $(1-\alpha)$-quantile of $\sqrt{ 2 L(\thetavt) - 2 L(\thetavs
)}$ is defined as
\begin{eqnarray*}
\zz_{\alpha} &\eqdef& \inf \bigl\{\zz\geq0 \colon \P \bigl(L(\thetavt) -
L\bigl(\thetavs\bigr) > \zz^{2}/2 \bigr) \leq\alpha \bigr\}
\\
&=& \inf \bigl\{\zz\geq0 \colon \E\Ind \bigl\{L(\thetavt) - L\bigl(\thetavs\bigr) >
\zz^{2}/2 \bigr\} \leq \alpha \bigr\}.
\end{eqnarray*}
Introduce for $x\geq0$ and $z,\Delta>0$ the following function:
%
\begin{equation}
\label{def:gDelta_xz} g_{\Delta}(x,z)\eqdef g \biggl( \frac{1}{2\Delta z} \bigl(
x^{2}-z^{2} \bigr) \biggr),
\end{equation}
where $g(x)$ is a three times differentiable nonnegative function,
and grows monotonously from $0$ to $1$, $g(x)=0$ for $x\leq0$
and $g(x)=1$ for $x\geq1$, therefore,
\[
\Ind \{ x> 1 \} \leq g(x)\leq\Ind \{ x> 0 \}\leq g(x+1).
\]
An example of such function is given in (A.8) in \citet
{SpZh2014PMBsupp}. It holds
\[
\Ind\{x-z > \Delta\} \leq g_{\Delta}(x,z)\leq\Ind(x-z>0)\leq
g_{\Delta
}(x,z+\Delta).
\]
This approximation is used in the proofs of Theorems \ref{thm:cumulat},
\ref{thm:cumulat2} and \ref{thm:resbias} in the part of Gaussian
approximation of Euclidean norm of a sum of independent vectors [see
Section~A.2 in \citet{SpZh2014PMBsupp}] yielding the error
rate $(\dimp^{3}/n)^{1/8}$ in the final bound [Theorems \ref
{thm:cumulat}, \ref{thm:cumulat2} and A.1 in
\citet{SpZh2014PMBsupp}].
The next result shows that the use of a smoothed quantile function
helps to improve
the accuracy of bootstrap approximation: it becomes
$(\dimp^{3}/n)^{1/2}$ instead of $ (\dimp^{3}/n)^{1/8} $.
The reason is that we do not need to account for the error induced by a
smooth approximation of the indicator function.

\begin{theorem}[{[Validity of the bootstrap in the smoothed case under
$(\SmB)$ condition]}]
\label{thm:smooth}
Let the conditions of Section~\ref{sect_conditions} be fulfilled.
It holds
for $\zz\geq  \max\{2,\sqrt{\dimp}\}+\CONST(\dimp+\yy)/\sqrt{n}$ and
$\Delta\in(0,0.22]$ with probability $\geq1-12\ex^{-\xx}$:
\[
\bigl\llvert \E g_{\Delta} \bigl(\sqrt{ 2L(\thetavt)-2L\bigl(\thetavs\bigr) }, \zz
\bigr) - \Eb g_{\Delta} \bigl(\sqrt{ 2\Lb\bigl(\thetavbt\bigr)-2\Lb(\thetavt) }, \zz
\bigr) \bigr\rrvert \leq \Delta_{\sm},
\]
where $\Delta_{\sm}\leq\CONST\{(\dimp+\yy)^{3}/n\}^{1/2}\Delta^{-3}$
in the case \ref{typical_local}.
An explicit definition of $\Delta_{\sm}$ is given in \textup{(B.59)},
\textup{(B.60)} in Section~\textup{B.6} of the
supplementary material [\citet{SpZh2014PMBsupp}].
\end{theorem}

The modified bootstrap quantile function reads as
\[
\label{def:zzbalpha_smooth} \zzb_{\Delta, \alpha} \eqdef \min \bigl\{\zz\geq0 \colon \Eb
g_{\Delta} \bigl(\sqrt{2\Lb\bigl(\thetavbt\bigr) - 2\Lb(\thetavt)}, \zz \bigr) \leq
\alpha \bigr\}.
\]
%
\section{Numerical results}
\label{sect:numer}
This section illustrates the performance of the multiplier bootstrap
for some artificial examples.
We especially aim to address the issues of noise misspecification and
of increasing modelling bias. It should be mentioned that the obtained
results are nicely consistent with the theoretical
statements.

In all the experiments, we took $10^{4}$ data samples for estimation
of the empirical c.d.f. of $\sqrt{2L(\thetavt)-2L(\thetavs)}$, and $
10^{4}$ $\{u_{1},\ldots,u_{n}\}$ samples for each of the $10^{4}$
data samples for the estimation of the quantiles of $\sqrt{2\Lb
(\thetavbt)-2\Lb(\thetavt)}$.

\subsection{Computational error}
\label{sect:numer_comp_err}
Here, we check numerically how well the multiplier procedure works in
the case of the correct model. Here, the modelling bias term $\|\Hc
^{-1}\Bsmb^{2}\Hc^{-1}\|$ from the $(\SmB)$ condition equals to
zero by its definition.
Let the data come from the following model:
$Y_{i}= \Psi_{i}^{\T}\thetav_{0}+\varepsilon_{i}$, for $i=1, \ldots,
n$, 
where $\varepsilon_{i}\sim\mathcal{N}(0,1)$, $\Psi_{i}\eqdef
(1,X_{i},X_{i}^{2},\ldots,X_{i}^{\dimp-1} )^{\T}$, the design
points $X_{1},\ldots,X_{n}$ are equidistant on $[0,1]$, and the
parameter vector $\thetav_{0}=(1,\ldots,1)^{\T}\in\R^{\dimp}$.
The true likelihood function is
$L(\thetav)= -\sum_{i=1}^{n}(Y_{i}-\Psi_{i}^{\T}\thetav)^{2}/2$. In
this experiment, we consider three cases: the scalar parameter $\dimp
= 1$, and the multivariate parameter $\dimp=3, 10$.

Table~\ref{tab:iid_Stnorm} shows the effective coverage probabilities
of the quantiles estimated using the multiplier bootstrap. The second
line contains the range of the nominal confidence levels: $0.99,\ldots
,0.75$. The first left column shows the sample size $n$ and the
second column---the parameter's dimension $\dimp$. The third left
column describes the distribution of the bootstrap weights: $
2\operatorname{Bernoulli}(0.5)$, $\mathcal{N}(1,1)$ or $\exp(1)$.
Below its second line, the table contains the frequencies of the
event: ``the real likelihood ratio $\leq$ the quantile of the
bootstrap likelihood ratio.''
%
\begin{table}
\caption{Coverage probabilities for the correct model}
\label{tab:iid_Stnorm}
\begin{tabular*}{\textwidth}{@{\extracolsep{\fill}}lcccccccc@{}}
\hline
&&&\multicolumn{6}{c@{}}{\textbf{Confidence levels}}\\[-6pt]
&&&\multicolumn{6}{c@{}}{\hrulefill}\\
 $\bolds{n}$
& $\bolds{\dimp}$ &$\bolds{{\mathcal{L}(u_{i})}}$ &$\mathbf{0\bolds{.}99}$ &$\mathbf
{0\bolds{.}95}$ &$\mathbf{0\bolds{.}90}$& $\mathbf{0\bolds{.}85}$& $\mathbf{0\bolds{.}80}$& $
\mathbf{0\bolds{.}75}$\\
\hline
{$50$}&
\phantom{0}{$1$}
&$2\operatorname{Bernoulli}(0.5)$
& $0.986$& $0.942$ & $0.892$ &$0.838$ &$0.792$ & $0.745$
\\
&& $\mathcal{N}(1,1)$
& $0.988$& $0.945$ & $0.895$ &$0.847$ &$0.803$ & $0.751$
\\
&&$\operatorname{exp}(1) $
& $0.988$& $0.942$ & $0.885$ &$0.833$ &$0.784$ & $0.729$
\\ [3pt]
{$50$}
&\phantom{0}{$3$}
&$2\operatorname{Bernoulli}(0.5)$
& $0.984$& $0.938$ & $0.885$ &$0.838$ &$ 0.788$ & $0.736$
\\
&& $\mathcal{N}(1,1)$
& $0.994$& $0.949$ & $0.897$ &$0.844$ &$0.789$ & $0.736$
\\
&&$\operatorname{exp}(1) $
& $0.984$& $0.917$ & $0.835$ &$0.776$ &$0.707$ & $0.650$
\\[3pt]
{$50$}
&{$10$}
&$2\operatorname{Bernoulli}(0.5)$
& $0.975$& $0.923$ & $0.866$ &$0.813$ &$0.764$ & $0.715$
\\
&& $\mathcal{N}(1,1)$
& $0.996$& $0.950$ & $0.877$ &$0.780$ &$0.721$ & $0.644$
\\
&&$\operatorname{exp}(1) $
& $0.952$& $0.827$ & $0.710$ &$0.617$ &$0.541$ & $0.473$
\\
\hline
\end{tabular*}
\end{table}

\subsection{Linear regression with misspecified heteroscedastic errors}
Here, we show on a linear regression model that the quality of the
confidence sets obtained by the multiplier bootstrap procedure is not
significantly deteriorated by misspecified heteroscedastic errors. Let
the data be defined as $Y_{i}= \Psi_{i}^{\T}\thetav_{0}+\sigma_{i}
\varepsilon_{i}$, $i=1, \ldots, n$. The i.i.d. random variables $
\varepsilon_{i}\sim \operatorname{Laplace}(0,2^{-1/2})$ are s.t. $\E(\varepsilon
_{i})=0$, $\Var(\varepsilon_{i})=1$. The coefficients $\sigma_{i}$
are deterministic: $\sigma_{i}\eqdef0.5 \{4-i (\mathrm{mod}\ 4)\}$.
The regressors $\Psi_{i}$ are the same as in the experiment \ref
{sect:numer_comp_err}. The quasi-likelihood function is also the same
as in the previous section: $L(\thetav)=-\sum_{i=1}^{n}(Y_{i}-\Psi
_{i}^{\T}\thetav)^{2}/2$, and it is misspecified, since it corresponds
to $\sigma_{i} \varepsilon_{i}\sim\mathcal{N}(0,1)$. The target
point $\thetavs=\thetav_{0}$, therefore, the modelling bias term $\|
\Hc^{-1}\Bsmb^{2}\Hc^{-1}\|$ from the $(\SmB)$ condition equals
to zero.

Here, we also consider three different parameter's dimensions: $
\dimp=1,3,10$ with $\thetav_{0}=(1,\ldots,1)^{\T}\in\R^{\dimp}$.
Table~\ref{tab:het_Lap} describes the second experiment's results
similarly to the Table~\ref{tab:iid_Stnorm}.%
%
\begin{table}
\caption{Coverage probabilities for case of misspecified heteroscedastic noise}

\label{tab:het_Lap}
\begin{tabular*}{\textwidth}{@{\extracolsep{\fill}}lcccccccc@{}}
\hline
&&&\multicolumn{6}{c@{}}{\textbf{Confidence levels}}\\[-6pt]
&&&\multicolumn{6}{c@{}}{\hrulefill}\\
 $\bolds{n}$
& $\bolds{\dimp}$ &$\bolds{{\mathcal{L}(u_{i})}}$ &$\mathbf{0\bolds{.}99}$ &$\mathbf
{0\bolds{.}95}$ &$\mathbf{0\bolds{.}90}$& $\mathbf{0\bolds{.}85}$& $\mathbf{0\bolds{.}80}$& $
\mathbf{0\bolds{.}75}$\\
\hline
\phantom{0}{$50$}
&\phantom{0}{$1$}
&$2\operatorname{Bernoulli}(0.5)$
& $0.988$& $0.947$ & $0.896$ &$0.849$ &$0.799$ & $0.752$
\\
&& $\mathcal{N}(1,1)$
& $0.990$& $0.949$ & $0.893$ &$0.844$ &$0.794$ & $0.746$
\\
&&$\exp(1)$
& $0.989$& $0.941$ & $0.881$ &$0.825$ &$0.770$ & $0.714$\\[3pt]
\phantom{0}{$50$}
&\phantom{0}{$3$}&$2\operatorname{Bernoulli}(0.5)$
& $0.984$& $0.937$ & $0.885$ &$0.834$ &$0.788$ & $0.739$
\\
&&$\mathcal{N}(1,1)$
& $0.996$& $0.955$ & $0.897$ &$0.839$ &$0.780$ & $0.722$
\\
&& $\exp(1)$
& $0.988$& $0.924$ & $0.846$ &$0.765$ &$0.701$ & $0.634$
\\ [3pt]
\phantom{0}{$50$}
&{$10$}&$2\operatorname{Bernoulli}(0.5)$
& $0.976$& $0.927$ & $0.870$ &$0.815$ &$0.765$ & $0.715$
\\
&& $\mathcal{N}(1,1)$
& $0.998$& $0.959$ & $0.891$ &$0.810$ &$0.731$ & $0.655$
\\
&& $\exp(1)$
& $ 0.967$& $0.850$ & $0.726$ &$0.630$ &$0.552$ & $0.479$
\\[3pt]
{$100$}
&{$10$}&$2\operatorname{Bernoulli}(0.5)$
& $ 0.985$& $0.935$ & $0.885$ &$0.833$ &$0.781$ & $0.733$
\\
&& $\mathcal{N}(1,1)$
& $0.998$& $0.970$ & $0.917$ &$0.857$ &$0.786$ & $0.723$
\\
&& $\exp(1)$
& $0.989$& $0.921$ & $0.826$ &$0.741$ &$ 0.663$ & $0.591$
\\ \hline
\end{tabular*}
\end{table}

One can see from the Tables~\ref{tab:iid_Stnorm} and \ref{tab:het_Lap}
that the bootstrap procedure does a good job even for small or moderate
samples like
50 or 100 if the parameter dimension is not too large.
The results are stable w.r.t. the noise misspecification.

The Rademacher and Gaussian weights demonstrate nearly the same nice
performance while
the procedure with exponential weights tends to underestimate the real
quantiles.
This effect becomes especially prominent when the parameter dimension
grows to 10.

\subsection{Biased constant regression with misspecified errors}
\label{examp:bias}
\begin{figure}[b]
\centering
\begin{tabular}{@{}c@{\quad}c@{}}

\includegraphics{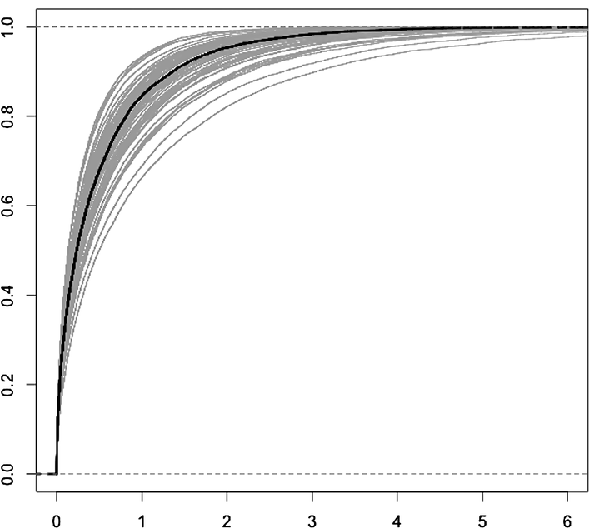}  & \includegraphics{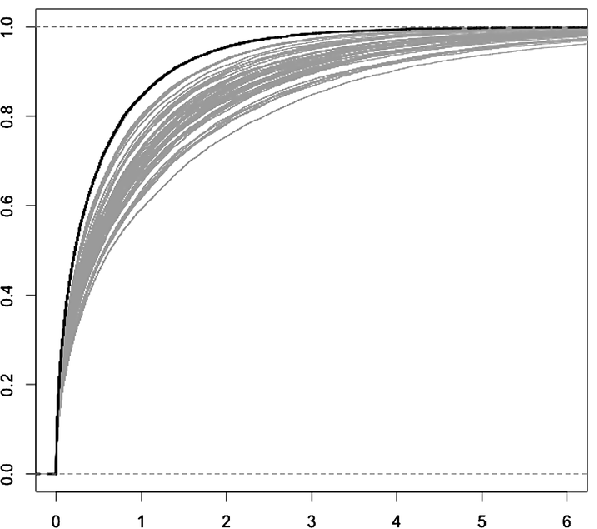}\\
\footnotesize{$Y_{i}=0.25\sin(X_{i})+\operatorname{Lap}(0,2^{-1/2}),  n=50$} &
\footnotesize{$Y_{i}=1.25\sin(X_{i})+\operatorname{Lap}(0,2^{-1/2}), n=50$}
\end{tabular}
\caption{Empirical distribution functions of the likelihood ratios.
$\protect\rule[2pt]{17pt}{0.5pt}$ Empirical distribution function of $L(\thetavt)-L(\thetavs)$ estimated
with $10^{4}~\Yv$ samples. \protect\includegraphics{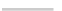}
50 empirical distribution functions of
$L^{\circ}(\tilde{\bolds{\theta}}^{\circ})-L^{\circ}(\tilde{\bolds{\theta}})$ estimated with $10^{4} \{u_i\}\sim\exp(1)$
samples.}\label{tab:biased}
\end{figure}

In the third experiment, we consider biased regression with
misspecified i.i.d. errors:
\begin{eqnarray*}
&&Y_{i}=\beta\sin(X_{i})+\varepsilon_{i},\qquad
\varepsilon_{i}\sim \operatorname{Laplace}\bigl(0, 2^{-1/2}\bigr),\qquad \mbox{i.i.d.},
\\
&& X_{i} \mbox{ are equidistant in } [0,2\pi].
\end{eqnarray*}
Taking the likelihood function $L(\thetav)=-\sum_{i=1}^{n}(Y_{i}-\thetav)^{2}/2$ yields $\thetavs=0$. Therefore, the
larger is the deterministic amplitude $\beta>0$, the bigger is bias
of the mean constant regression. The $(\SmB)$ condition reads as
\begin{eqnarray*}
\bigl\|\Hc^{-1}\Bsmb^{2}\Hc^{-1}\bigr\|&=&1-
\frac{\sum_{i=1}^{n} \Var
Y_{i}}{\beta^{2}\sum_{i=1}^{n}\sin^{2}(X_{i})+\sum_{i=1}^{n}\Var Y_{i}} %
\\
&=& 1-\frac{1}{\beta^{2}(n-1)/2n+1}
\\
&\leq& 1/\sqrt{n}.
\end{eqnarray*}
Consider the sample size $n=50$, and two cases: $\beta=0.25$ with
fulfilled $(\SmB)$ condition and $\beta=1.25$ when $(\SmB)$ does not hold. Table~\ref{tab:sinnLap} shows that for the
large bias quantiles yielded by the multiplier bootstrap are conservative.
This conservative property of the multiplier bootstrap quantiles is
also illustrated with the graphs in Figure~\ref{tab:biased}.
They show the empirical distribution functions of the likelihood ratio
statistics $L(\thetavt)-L(\thetavs)$ and $\Lb(\thetavbt)-\Lb(\thetavt
)$ for $\beta=0.25$ and $\beta=1.25$. On the right graph for $
\beta=1.25$ the empirical distribution functions for the bootstrap
case are smaller than the one for the $\Yv$ case. It means that for
the large bias the bootstrap quantiles are bigger than the $\Yv$
quantiles, which increases the diameter of the confidence set based on
the bootstrap quantiles. This confidence set remains valid, since it
still contains the true parameter with a given confidence level.

%
\begin{table}
\caption{Coverage probabilities for the noise-misspecified biased regression}
\label{tab:sinnLap}
\begin{tabular*}{\textwidth}{@{\extracolsep{\fill}}lcccccccc@{}}
\hline
&&&\multicolumn{6}{c@{}}{\textbf{Confidence levels}}\\[-6pt]
&&&\multicolumn{6}{c@{}}{\hrulefill}\\
$ \bolds{n} $&
$\bolds{\mathcal{L}(u_{i})}$& $\bolds{\beta}$&$\mathbf{0\bolds{.}99}$ &$\mathbf{0\bolds{.}95}$
&$\mathbf{0\bolds{.}90}$& $\mathbf{0\bolds{.}85}$& $\mathbf{0\bolds{.}80}$& $\mathbf
{0\bolds{.}75}$\\
\hline
$50$&
$\mathcal{N}(1,1)$
& $0.25$&
$0.98$ & $0.94$& $0.89$ & $0.84$ &$0.79$ &$0.74$
\\
&&$1.25$&$1.0$\phantom{0}&$0.99$&$0.97$&$0.94$&$0.91$&$0.87$\\
\hline
\end{tabular*}
\end{table}

%

%

%
%

\begin{figure}

\includegraphics{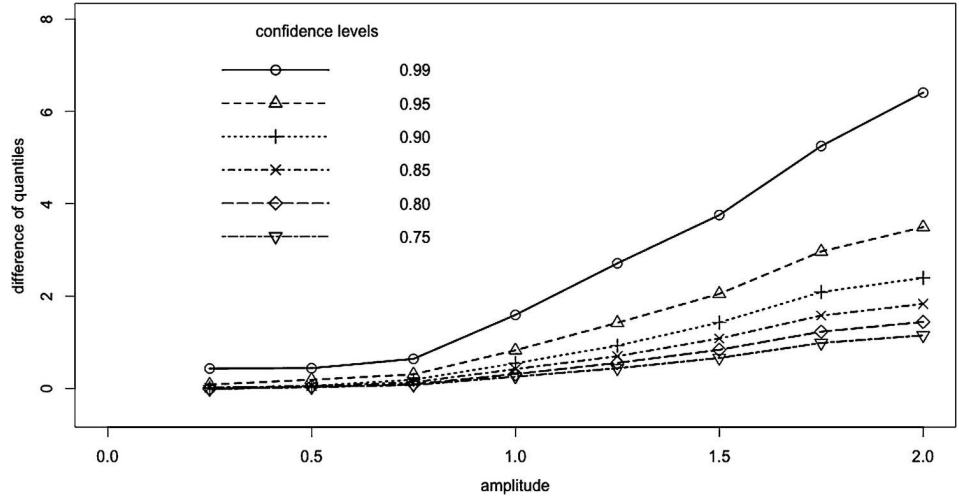}

\caption{The difference  (``Bootstrap
quantile''$-$``$\Yv$-quantile'') growing with modelling bias.}
\label{fig:bias}
\end{figure}

Figure~\ref{fig:bias} shows the growth
of the difference between the quantiles of $\Lb(\thetavbt)-\Lb(\thetavt
)$ and $L(\thetavt)-L(\thetavs)$ with increasing $\beta$ for the
range of the confidence levels: $0.75, 0.8, \ldots, 0.99$.
%

\subsection{Logistic regression with bias}
\label{examp:logistic}
\begin{figure}
\centering
\begin{tabular}{@{}cc@{}}

\includegraphics{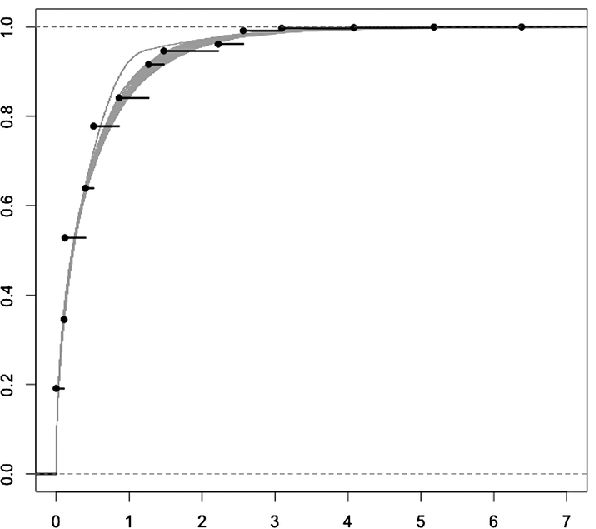}  & \includegraphics{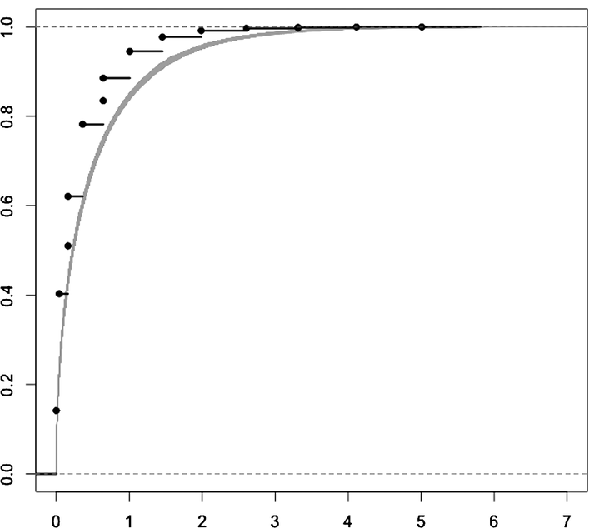}\\
\footnotesize{$Y_{i}\sim \operatorname{Bernoulli}(0.1 X_{i}), n=50$} & \footnotesize{$Y_{i}\sim \operatorname{Bernoulli}(0.5 X_{i}), n=50$}
\end{tabular}
\caption{Empirical distribution functions of the likelihood ratios for logistic regression. $\protect\rule[2pt]{17pt}{0.5pt}$ Empirical distribution function of $L(\thetavt)-L(\thetavs)$ estimated
with $10^{4}~\Yv$ samples. \protect\includegraphics{1355i01.eps}
50 empirical distribution functions of
$L^{\circ}(\tilde{\bolds{\theta}}^{\circ})-L^{\circ}(\tilde{\bolds{\theta}})$ estimated with $10^{4} \{u_i\}\sim\exp(1)$
samples.}\label{tab:logistic}
\end{figure}

In this example, we consider logistic regression. Let the data come
from the following distribution:
\[
Y_{i}\sim \operatorname{Bernoulli}(\beta X_{i}),\qquad X_{i}
\mbox{ are equidistant in }[0,2], \beta\in(0,1/2].
\]
Consider the likelihood function corresponding to the i.i.d. observations
\[
L(\thetav)=\sum_{i=1}^{n} \bigl
\{Y_{i}\thetav-\log\bigl(1+\ex ^{\thetav}\bigr) \bigr\}.
\]
By definition \eqref{thetavs} $\thetavs=\log\{\beta/(1-\beta)\}$,
bigger values of $\beta$ induce larger modelling bias. Indeed, the
$(\SmB)$ condition reads as
\begin{eqnarray*}
\bigl\|\Hc^{-1}\Bsmb^{2}\Hc^{-1}\bigr\|&=&
\frac{\beta^{2}\sum_{i=1}^{n}(X_{i}-1)^{2}}{n\beta^{2}+\beta(1-2\beta)\sum_{i=1}^{n}X_{i}}
\\
&=&\frac{\beta}{1-\beta}\cdot\frac{n+1}{3(n-1)}
\\
&\leq& 1/\sqrt{n}.
\end{eqnarray*}
The graphs on Figure~\ref{tab:logistic} demonstrate the
conservativeness of bootstrap quantiles. Here, we consider two cases: $
\beta=0.1$ and $\beta=0.5$. Similarly to the Example~\ref
{examp:bias} in the case of the bigger $\beta$ on the right graph of
Figure~\ref{tab:logistic}, the empirical distribution functions of $\Lb
(\thetavbt)-\Lb(\thetavt)$ are smaller than the one for $L(\thetavt
)-L(\thetavs)$.
%

%
%
%
%

\section{Conditions}
\label{sect_conditions}
Here, we state the conditions required for the main results. The
conditions in Section~\ref{sect:ConditGeneral} come from the general
finite sample theory by \citet{Sp2012Pa}. They are required for the
results of Sections~B.1 and B.2 in the
supplementary material [\citet{SpZh2014PMBsupp}].
The conditions in Section~\ref{sect:ConditAddBoot} are necessary to
prove the results on multiplier bootstrap from Section~\ref{sect_mainres}.
In Section~B.3 in \citet{SpZh2014PMBsupp}, we
consider these conditions in detail for several examples: i.i.d.
observations, generalized linear model and linear quantile regression.

\subsection{Basic conditions}
\label{sect:ConditGeneral}
Introduce the stochastic part of the likelihood process:
$\zeta(\thetav)\eqdef L(\thetav)- \E L(\thetav)$, and its marginal
summand: $\zeta_{i}(\thetav)\eqdef\ell_{i}(\thetav)-\E\ell_{i}(\thetav
)$.
\begin{longlist}
\item[$(\ED_{0})$]
There exist a positive-definite symmetric matrix $\VPc^{2}$ and
constants $\gm>0, \nunu\geq1$ such that $\Var \{\nabla_{\thetav
}\zeta(\thetavs) \}\leq\VPc^{2} $ and
\[
\sup_{\gammav\in\R^{\dimp}} \log\E\exp \biggl\{ \lambda\frac{\gammav^{\T} \nabla_{\thetav}\zeta\bigl(\thetavs\bigr)}{\|\VPc\gammav
\|}
\biggr\}\leq\nunu^{2}\lambda^{2}/2, \qquad |\lambda|\leq\gm.
\]
\end{longlist}
\begin{longlist}
\item[$(\ED_{2})$]
There exist a constant $\omega> 0$ and for each $\rr>0$ a
constant $\gm_{2}(\rr)$ such that it holds for all $\thetav\in
\Thetas(\rr)$ and for $j=1,2$
\[
\sup_{{\gammav_{j}\in\R^{\dimp}\\
\|\gammav_{j}\|\leq1}
} \log\E\exp \biggl\{ \frac{\lambda}{\omega}
\gammav_{1}^{\T}\DPc^{-1}\nabla_{\thetav}^{2}
\zeta (\thetav)\DPc^{-1}\gammav_{2} \biggr\} \leq
\nunu^{2}\lambda^{2}/2,\qquad |\lambda|\leq\gm_{2}(
\rr).
\]
\end{longlist}
\begin{longlist}
\item[$(\LL_{0})$]
For each $\rr\in[0,\rups]$ [$\rups$ comes from condition
(B.1) of Theorem~B.1 in \citet
{SpZh2014PMBsupp}] there exists a constant $\rddelta(\rr)\in[0,1/2]$
s.t. for all $\thetav\in\Thetas(\rr)$ it holds
\[
\bigl\|\DPc^{-1}D^{2}(\thetav)\DPc^{-1}-
\Id_{\dimp}\bigr\|\leq\rddelta(\rr),
\]
where $D^{2}(\thetav)\eqdef-\nabla_{\thetav}^{2}\E L(\thetav)$, $\Thetas(\rr)\eqdef \{\thetav: \|\DPc(\thetav-\thetavs)\|\leq\rr
 \}$.
\end{longlist}
\begin{longlist}
\item[$(\mathcal{I})$]
There exists a constant $\gmu>0$ s.t.
$\gmu^{2} \DPc^{2}\geq\VPc^{2}$.
\end{longlist}
\begin{longlist}
\item[$(\LL\rr)$]
For each $\rr>\rups$ there exists a value $\gmi(\rr)>0$ s.t.
$\rr\gmi(\rr)\rightarrow+\infty$ for $\rr\rightarrow+\infty$ and
$\forall\thetav: \|\DPc(\thetav-\thetavs)\|=\rr$ it holds
\[
-2 \bigl\{\E L(\thetav)-\E L\bigl(\thetavs\bigr) \bigr\}\geq\rr^{2}\gmi(\rr).
\]
\end{longlist}
%
\subsection{Conditions required for the bootstrap validity}
\label{sect:ConditAddBoot}
\begin{longlist}
\item[$(\SmB)$] 
There exists a constant $\deltasmb^{2}\in[0,1/8]$
such that it holds for the matrices $\Hc^{2}$, $\Bsmb^{2}$ defined
in \eqref{def:Hc} and \eqref{def:Bsmb}.
\[
\bigl\|\Hc^{-1}\Bsmb^{2}\Hc^{-1}\bigr\|\leq
\deltasmb^{2}\leq\CONST\dimp n^{-1/2}.
\]
\end{longlist}
\begin{longlist}
\item[$(\ED_{2m})$]
For each $\rr>0$, $i=1,\ldots,n$, $j=1,2$ and for all $
\thetav\in\Thetas(\rr)$ it holds for the values $\omega\geq0$ and
$\gm_{2}(\rr)$ from the condition $(\ED_{2})$
\[
\mathop{\mathop{\sup_{{\gammav_{j}\in\R^{\dimp}}}}_{
\|\gammav_{j}\|\leq1}
} \log\E\exp \biggl\{ \frac{\lambda}{\omega}
\gammav_{1}^{\T}\DPc^{-1} \nabla_{\thetav}^{2}
\zeta_{i}(\thetav) \DPc^{-1}\gammav_{2} \biggr\}
\leq\frac{\nunu^{2}\lambda^{2}}{2n},\qquad |\lambda|\leq\gm _{2}(\rr).
\]
\end{longlist}
\begin{longlist}
\item[$(\LL_{0m})$] 
For each $\rr>0$, $i=1,\ldots,n$ and for all $\thetav\in
\Thetas(\rr)$, there exists a constant $\CONST_{m}(\rr)\geq0$ such that
\[
\bigl\|\DPc^{-1}\nabla_{\thetav}^{2}\E\ell_{i}(
\thetav)\DPc^{-1}\bigr\|\leq {\CONST_{m}(\rr)} {n}^{-1}.
\]
\end{longlist}
\begin{longlist}
\item[$(\mathcal{I}_{B})$] 
There exists a constant $\gmu_{B}^{2}\geq0$ s.t.
$\gmu_{B}^{2} \DPc^{2}\geq\Bsmb^{2}$.
\end{longlist}
\begin{longlist}
\item[$(\SD_{1})$] 
There exists a constant $0\leq\delta_{v}^{2}\leq\CONST\dimp
/n$. such that it holds for all $i=1,\ldots,n$ with exponentially
high probability
\[
\bigl\llVert \Hc^{-1} \bigl\{\nabla_{\thetav}
\ell_{i}\bigl(\thetavs\bigr)\nabla_{\thetav
}\ell_{i}\bigl(
\thetavs\bigr)^{\T}- \E \bigl[\nabla_{\thetav}\ell_{i}\bigl(
\thetavs \bigr)\nabla_{\thetav}\ell_{i}\bigl(\thetavs\bigr)^{\T} \bigr]
\bigr\}\Hc^{-1}\bigr\rrVert \leq\delta_{v}^{2}.%
\]
\end{longlist}
\begin{longlist}
\item[$(Eb)$] 
The bootstrap weights $u_{i}$ are i.i.d., independent of the
data $\Yv$, and
\begin{eqnarray*}
&&\E u_{i} =1, \qquad\Var u_{i}=1,
\\
&&\log\E\exp \bigl\{ \lambda(u_{i}-1) \bigr\}\leq\nunu^{2}
\lambda^{2}/2,\qquad |\lambda|\leq\gm.
\end{eqnarray*}
\end{longlist}

\subsection{Small modelling bias condition for some models}
\label{sect:smb_condit}
Here, we specify the condition $(\SmB)$ for some particular
models. If the observations $Y_{1},\ldots, Y_{n}$ are i.i.d., then $
\nabla_{\thetav}\E L(\thetavs)=n\nabla_{\thetav}\E\ell_{i}(\thetavs
)=0$, and $\Bsmb^{2}=0$. The next example is the generalized linear
model: the parametric probability distribution family $ \{\P
_{\upsilon} \}$ is an exponential family with a canonical
parameterization. The log-density for this family can be expressed as
\[
\ell(\upsilon)=yv-\dv(\upsilon)
\]
for a convex function $\dv(\cdot)$. Table~\ref{tab:glm_main_cases}
provides some examples of $ \{\P_{\upsilon} \}$ and $\dv
(\cdot)$.
Taking $ \{\P_{\upsilon} \}$ as a parametric family and $
\Psi_{i}^{\T}\thetav$ as linear predictors for some deterministic
regressors $\Psi_{i}\in\R^{\dimp}$ yields the following quasi
log-likelihood function:
\[
L(\thetav)=\sum_{i=1}^{n} \bigl
\{Y_{i}\Psi_{i}^{\T}\thetav-\dv \bigl(
\Psi_{i}^{\T}\thetav\bigr) \bigr\}.
\]
%
\begin{table}
\caption{Examples of the GLM}
\label{tab:glm_main_cases}
\begin{tabular*}{\textwidth}{@{\extracolsep{\fill}}lcc@{}}
\hline
$\bolds{\P_{\upsilon}}$ &$\bolds{\dv(\upsilon)}$ & $\bolds{\dv^{\prime}(\upsilon)}$ \textbf{(natural parameter)} \\
\hline
${\mathcal{N}(\upsilon,1)}$
& ${\upsilon^{2}/2}$& ${\upsilon}$
\\
$\operatorname{Exp}(-\upsilon)$
& $-\log(-\upsilon)$& $-1/\upsilon$
\\
$\operatorname{Pois}(\ex^{\upsilon})$
& $\ex^{\upsilon}$& $\ex^{\upsilon}$
\\
$\operatorname{Binom}(1,\frac{\ex^{\upsilon}}{\ex^{\upsilon}+1} )$
& $\log (\ex^{\upsilon}+1 )$& $\frac{\ex^{\upsilon}}{\ex
^{\upsilon}+1}$ \\
\hline
\end{tabular*}
\end{table}
\noindent It holds
\begin{eqnarray*}
&& \bigl\|\Hc^{-1}\Bsmb^{2}\Hc^{-1}\bigr\|
\\
&&\qquad\leq 1- \min_{1\leq i\leq n} \frac{\Var Y_{i}}{\Var Y_{i} + \{\E Y_{i}-\dv^{\prime}(\Psi_{i}^{\T
}\thetavs) \}^{2}} \in[0,1).
\end{eqnarray*}
It is important that $\E_{\thetavs}Y_{i}=\dv^{\prime}(\Psi_{i}^{\T
}\thetavs)$, that is, in the case of the correct parametric model $\P
\in \{\P_{\upsilon} \}$ the modelling bias is indeed equal to zero.

Now let us consider the linear quantile regression. Let the
observations $Y_{1},\ldots,Y_{n}$ be scalar, and the design points $
X_{1},\ldots,X_{n}$ be deterministic. Let $\tau\in(0,1)$ denote a
fixed known quantile level. The object of estimation is a quantile
function $q_{\tau}(x)$ s.t.
\[
\P \bigl(Y_{i}< q_{\tau}(X_{i}) \bigr) =\tau\qquad
\forall i=1,\ldots,n.
\]
Using the quantile regression approach by \citet{Koenker1978regression},
this problem can be treated with the quasi maximum likelihood method
and the following log-likelihood function:
%
\begin{eqnarray}
\label{def:qr} %
L(\thetav)&=& -\sum_{i=1}^{n}
\rho_{\tau} \bigl(Y_{i}-\Psi_{i}^{\T
}
\thetav \bigr),
\nonumber
\\[-8pt]
\\[-8pt]
\nonumber
\rho_{\tau}(x)&\eqdef& x \bigl(\tau-\Ind \{x<0 \} \bigr), %
\end{eqnarray}
where $\Psi_{i}\in\R^{\dimp}$ are known regressors. This
log-likelihood function corresponds to asymmetric Laplace distribution
with the density $\tau(1-\tau)\ex^{-\rho_{\tau}(x-a)}$.
It holds
\begin{eqnarray*}
&& \bigl\|\Hc^{-1}\Bsmb^{2}\Hc^{-1}\bigr\|
\\
&&\qquad\leq 1- \min_{1\leq i\leq n} \frac{\Var (\tau-\Ind \{Y_{i}-\Psi_{i}^{\T}\thetavs<0 \}
 )}{\Var (\tau-\Ind \{Y_{i}-\Psi_{i}^{\T}\thetavs<0
\} )+ (\tau-\P \{Y_{i}-\Psi_{i}^{\T}\thetavs<0 \} )^{2}}.
\end{eqnarray*}
If $\P \{Y_{i}-\Psi_{i}^{\T}\thetavs<0 \}\equiv\tau$, then
the right-hand side of the last inequality is equal to zero.

\subsection{Dependence of the involved terms on the sample size and
parameter dimension}
\label{typical_local}
Here, we consider the case of the i.i.d. observations $Y_{1},\ldots,
Y_{n}$ and $\yy= \CONST\log{n}$ in order to specify the dependence
of the nonasymptotic bounds on $n$ and $\dimp$.
In Section~B.3 of the supplementary material [\citet
{SpZh2014PMBsupp}], we also consider generalized linear model and
quantile regression.
Example~5.1 in \citet{Sp2012Pa} demonstrates that in this situation $\gm
=\CONST\sqrt{n}$ and $\omega=\CONST/\sqrt{n}$. This yields $\ZZ(\yy
)= \CONST\sqrt{\dimp+\yy}$ for some constant $\CONST\geq1.85$, for
the function $\ZZ(\yy)$ given in (B.4) in Section~B.1
of the supplementary material [\citet{SpZh2014PMBsupp}].
Similarly, it can be checked that $\gm_{2}(\rr)$ from condition $(\ED_{2})$
is proportional to $\sqrt{n}$: due to independence of the observations
\begin{eqnarray*}
&& \log\E\exp \biggl\{ \frac{\lambda}{\omega}\gammav_{1}^{\T}
\DPc^{-1}\nabla_{\thetav}^{2}\zeta (\thetav)
\DPc^{-1}\gammav_{2} \biggr\}
\\
&&\qquad=\sum_{i=1}^{n}\log\E\exp \biggl\{
\frac{\lambda}{\sqrt{n}}\frac{1}{\omega\sqrt{n}} \gammav_{1}^{\T}
\dpc^{-1}\nabla_{\thetav}^{2}\zeta_{i}(
\thetav)\dpc ^{-1}\gammav_{2} \biggr\}
\\
&&\qquad\leq n\frac{\lambda^{2}}{n}\CONST\qquad \mbox{for } |\lambda|\leq\bar {
\gm}_{2}(\rr)\sqrt{n},
\end{eqnarray*}
where $\zeta_{i}(\thetav)\eqdef\ell_{i}(\thetav)-\E\ell_{i}(\thetav)$, $\dpc^{2}\eqdef-\nabla^{2}_{\thetav}\E\ell_{i}(\thetavs)$ and $
\DPc^{2}=n\dpc^{2}$ in the i.i.d. case. Function $\bar{\gm}_{2}(\rr
)$ denotes the marginal analog of $\gm_{2}(\rr)$.

Let us show that for the value $\delta(\rr)$ from condition $(\LL_{0})$ it holds $\delta(\rr)=\CONST\rr/\sqrt{n}$. Suppose for all
$\thetav\in\Thetas(\rr)$ and $\gammav\in\R^{\dimp}:\|\gammav\|=1$
$\|\DPc^{-1}\gammav^{\T}\nabla_{\thetav}^{3}\E L(\thetav)\DPc^{-1}\|
\leq\CONST$, then it holds for some $\bar{\thetav}\in\Thetas(\rr)$
\begin{eqnarray*}
\bigl \|\DPc^{-1}D^{2}(\thetav)\DPc^{-1}-
\Id_{\dimp}\bigr\|&= &\bigl\|\DPc^{-1}\bigl(\thetavs-\thetav\bigr)^{\T}
\nabla_{\thetav}^{3}\E L(\bar{\thetav })\DPc^{-1}\bigr\|
\\
&=&\bigl\|\DPc^{-1}\bigl(\thetavs-\thetav\bigr)^{\T}\DPc\DPc^{-1}
\nabla_{\thetav
}^{3}\E L(\bar{\thetav})\DPc^{-1}\bigr\|
\\
&\leq& \rr\bigl\|\DPc^{-1}\bigr\|\bigl\|\DPc^{-1}\gammav^{\T}
\nabla_{\thetav}^{3}\E L(\bar{\thetav})\DPc^{-1}\bigr\|\leq
\CONST\rr/\sqrt{n}.
\end{eqnarray*}
Similarly, ${\CONST_{m}(\rr)}\leq\CONST\rr/\sqrt{n}+\CONST$ in
condition $(\LL_{0m})$.

The next remark helps to check the global identifiability condition $(\LL\rr)$ in many situations.
Suppose that the parameter domain $\Theta$ is compact and $n$ is
sufficiently large, then the value $\gmi(\rr)$ from condition $(\LL\rr)$ can be taken as $\CONST\{1-\rr/\sqrt{n}\}\approx\CONST$.
Indeed, for $\thetav: \|\DPc(\thetav-\thetavs)\|=\rr$
\begin{eqnarray*}
-2 \bigl\{\E L(\thetav)-\E L\bigl(\thetavs\bigr) \bigr\} &\geq& \rr^{2} \bigl
\{1-\rr\bigl\|\DPc^{-1}\bigr\|\bigl\|\DPc^{-1}\gammav^{\T}\nabla
_{\thetav}^{3}\E L(\bar{\thetav})\DPc^{-1}\bigr\| \bigr\}
\\
&\geq& \rr^{2}(1- \CONST\rr/\sqrt{n}).
\end{eqnarray*}
Due to the obtained orders, conditions (B.1) and
(B.19) of Theorems B.1 and B.6 (in the supplementary material [\citet{SpZh2014PMBsupp}])
on concentration of the MLEs $\thetavt,  \thetavbt$ require $\rups
\geq\CONST\sqrt{\dimp+\yy}$.


\section*{Acknowledgments}
The authors are very grateful to the anonymous
referees for their careful reading of the manuscript and many helpful
remarks and suggestions.

\begin{supplement}[id=suppA]
\stitle{Supplement to ``Bootstrap confidence sets under model
misspecification''}
\slink[doi]{10.1214/15-AOS1355SUPP} 
\sdatatype{.pdf}
\sfilename{aos1355\_supp.pdf}
\sdescription{The supplementary material contains
a proof of the square-root Wilks approximation for the bootstrap world,
proofs of the main results from Section~\ref{sect_mainres}, and results
on Gaussian approximation for $\ell_{2}$-norm of a sum of independent
vectors.}
\end{supplement}


\begin{thebibliography}{29}
\bibitem[\protect\citeauthoryear{Aerts and Claeskens}{2001}]{Aerts2001missp}
\begin{barticle}[mr]
\bauthor{\bsnm{Aerts},~\bfnm{Marc}\binits{M.}} \AND
\bauthor{\bsnm{Claeskens},~\bfnm{Gerda}\binits{G.}}
(\byear{2001}).
\btitle{Bootstrap tests for misspecified models, with application to clustered binary data}.
\bjournal{Comput. Statist. Data Anal.}
\bvolume{36}
\bpages{383--401}.
\bid{doi={10.1016/S0167-9473(00)00051-7}, issn={0167-9473}, mr={1843893}}
\end{barticle}
%
\bptok{imsref}%
\endbibitem\vspace*{-7pt}

\bibitem[\protect\citeauthoryear{Arlot, Blanchard and Roquain}{2010}]{ArlotBlanch2010}
\begin{barticle}[mr]
\bauthor{\bsnm{Arlot},~\bfnm{Sylvain}\binits{S.}},
\bauthor{\bsnm{Blanchard},~\bfnm{Gilles}\binits{G.}} \AND
\bauthor{\bsnm{Roquain},~\bfnm{Etienne}\binits{E.}}
(\byear{2010}).
\btitle{Some nonasymptotic results on resampling in high dimension. I. {C}onfidence regions}.
\bjournal{Ann. Statist.}
\bvolume{38}
\bpages{51--82}.
\bid{doi={10.1214/08-AOS667}, issn={0090-5364}, mr={2589316}}
\end{barticle}
%
\bptok{imsref}%
\endbibitem

\bibitem[\protect\citeauthoryear{Barbe and Bertail}{1995}]{Barbe1995weighted}
\begin{bbook}[mr]
\bauthor{\bsnm{Barbe},~\bfnm{Philippe}\binits{P.}} \AND
\bauthor{\bsnm{Bertail},~\bfnm{Patrice}\binits{P.}}
(\byear{1995}).
\btitle{The Weighted Bootstrap}.
\bseries{Lecture Notes in Statistics}
\bvolume{98}.
\bpublisher{Springer},
\blocation{New York}.
\bid{doi={10.1007/978-1-4612-2532-4}, mr={2195545}}
\end{bbook}
%
\bptok{imsref}%
\endbibitem

\bibitem[\protect\citeauthoryear{B{\"u}cher and Dette}{2013}]{Bucher2013multiplier}
\begin{barticle}[mr]
\bauthor{\bsnm{B{\"u}cher},~\bfnm{Axel}\binits{A.}} \AND
\bauthor{\bsnm{Dette},~\bfnm{Holger}\binits{H.}}
(\byear{2013}).
\btitle{Multiplier bootstrap of tail copulas with applications}.
\bjournal{Bernoulli}
\bvolume{19}
\bpages{1655--1687}.
\bid{doi={10.3150/12-BEJ425}, issn={1350-7265}, mr={3129029}}
\end{barticle}
%
\bptok{imsref}%
\endbibitem

\bibitem[\protect\citeauthoryear{Chatterjee and Bose}{2005}]{ChattBose2005generalized}
\begin{barticle}[mr]
\bauthor{\bsnm{Chatterjee},~\bfnm{Snigdhansu}\binits{S.}} \AND
\bauthor{\bsnm{Bose},~\bfnm{Arup}\binits{A.}}
(\byear{2005}).
\btitle{Generalized bootstrap for estimating equations}.
\bjournal{Ann. Statist.}
\bvolume{33}
\bpages{414--436}.
\bid{doi={10.1214/009053604000000904}, issn={0090-5364}, mr={2157808}}
\end{barticle}
%
\bptok{imsref}%
\endbibitem

\bibitem[\protect\citeauthoryear{Chen and Pouzo}{2009}]{Chen2009efficient}
\begin{barticle}[mr]
\bauthor{\bsnm{Chen},~\bfnm{Xiaohong}\binits{X.}} \AND
\bauthor{\bsnm{Pouzo},~\bfnm{Demian}\binits{D.}}
(\byear{2009}).
\btitle{Efficient estimation of semiparametric conditional moment models with possibly nonsmooth residuals}.
\bjournal{J. Econometrics}
\bvolume{152}
\bpages{46--60}.
\bid{doi={10.1016/j.jeconom.2009.02.002}, issn={0304-4076}, mr={2562763}}
\end{barticle}
%
\bptok{imsref}%
\endbibitem

\bibitem[\protect\citeauthoryear{Chen and Pouzo}{2015}]{Chen2014sieve}
\begin{barticle}[auto:parserefs-M02]
\bauthor{\bsnm{Chen},~\bfnm{X.}\binits{X.}} \AND
\bauthor{\bsnm{Pouzo},~\bfnm{D.}\binits{D.}}
(\byear{2015}).
\btitle{Sieve {W}ald and QLR inferences on semi/nonparametric conditional moment models}.
\bjournal{Econometrica}
\bvolume{83}
\bpages{1013--1079}.
\bid{mr={3357484}}
\end{barticle}
%
\bptok{imsref}%
\endbibitem

\bibitem[\protect\citeauthoryear{Chernozhukov, Chetverikov and Kato}{2013}]{ChernoMultBoot}
\begin{barticle}[mr]
\bauthor{\bsnm{Chernozhukov},~\bfnm{Victor}\binits{V.}},
\bauthor{\bsnm{Chetverikov},~\bfnm{Denis}\binits{D.}} \AND
\bauthor{\bsnm{Kato},~\bfnm{Kengo}\binits{K.}}
(\byear{2013}).
\btitle{Gaussian approximations and multiplier bootstrap for maxima of sums of high-dimensional random vectors}.
\bjournal{Ann. Statist.}
\bvolume{41}
\bpages{2786--2819}.
\bid{doi={10.1214/13-AOS1161}, issn={0090-5364}, mr={3161448}}
\end{barticle}
%
\bptok{imsref}%
\endbibitem

\bibitem[\protect\citeauthoryear{Efron}{1979}]{Efron1979}
\begin{barticle}[mr]
\bauthor{\bsnm{Efron},~\bfnm{B.}\binits{B.}}
(\byear{1979}).
\btitle{Bootstrap methods: Another look at the jackknife}.
\bjournal{Ann. Statist.}
\bvolume{7}
\bpages{1--26}.
\bid{issn={0090-5364}, mr={0515681}}
\end{barticle}
%
\bptok{imsref}%
\endbibitem

\bibitem[\protect\citeauthoryear{Hall}{1992}]{Hall1992bootstbook}
\begin{bbook}[mr]
\bauthor{\bsnm{Hall},~\bfnm{Peter}\binits{P.}}
(\byear{1992}).
\btitle{The Bootstrap and {E}dgeworth Expansion}.
\bpublisher{Springer},
\blocation{New York}.
\bid{doi={10.1007/978-1-4612-4384-7}, mr={1145237}}
\end{bbook}
%
\bptok{imsref}%
\endbibitem

\bibitem[\protect\citeauthoryear{Hall}{2005}]{Hall2005GMMbook}
\begin{bbook}[mr]
\bauthor{\bsnm{Hall},~\bfnm{Alastair~R.}\binits{A.~R.}}
(\byear{2005}).
\btitle{Generalized Method of Moments}.
\bpublisher{Oxford Univ. Press},
\blocation{Oxford}.
\bid{mr={2135106}}
\end{bbook}
%
\bptok{imsref}%
\endbibitem

\bibitem[\protect\citeauthoryear{Horowitz}{2001}]{Horowitz2001bootstrapHandbook}
\begin{barticle}[auto:parserefs-M02]
\bauthor{\bsnm{Horowitz},~\bfnm{J.~L.}\binits{J.~L.}}
(\byear{2001}).
\btitle{The bootstrap}.
\bjournal{Handbook of Econometrics}
\bvolume{5}
\bpages{3159--3228}.
\end{barticle}
%
\bptok{imsref}%
\endbibitem

\bibitem[\protect\citeauthoryear{Huber}{1967}]{huber1967}
\begin{bincollection}[mr]
\bauthor{\bsnm{Huber},~\bfnm{Peter~J.}\binits{P.~J.}}
(\byear{1967}).
\btitle{The behavior of maximum likelihood estimates under nonstandard conditions}.
In \bbooktitle{Proc. {F}ifth {B}erkeley {S}ympos. {M}ath. {S}tatist. and {P}robability ({B}erkeley, {C}alif., 1965/66), {V}ol. I: {S}tatistics}
\bpages{221--233}.
\bpublisher{Univ. California Press},
\blocation{Berkeley, CA}.
\bid{mr={0216620}}
\end{bincollection}
%
\bptok{imsref}%
\endbibitem

\bibitem[\protect\citeauthoryear{Janssen}{1994}]{Janssen1994weighted}
\begin{barticle}[mr]
\bauthor{\bsnm{Janssen},~\bfnm{Paul}\binits{P.}}
(\byear{1994}).
\btitle{Weighted bootstrapping of {$U$}-statistics}.
\bjournal{J. Statist. Plann. Inference}
\bvolume{38}
\bpages{31--41}.
\bid{doi={10.1016/0378-3758(92)00156-X}, issn={0378-3758}, mr={1256846}}
\end{barticle}
%
\bptok{imsref}%
\endbibitem

\bibitem[\protect\citeauthoryear{Janssen and Pauls}{2003}]{Janssen2003bootstrap}
\begin{barticle}[mr]
\bauthor{\bsnm{Janssen},~\bfnm{Arnold}\binits{A.}} \AND
\bauthor{\bsnm{Pauls},~\bfnm{Thorsten}\binits{T.}}
(\byear{2003}).
\btitle{How do bootstrap and permutation tests work?}
\bjournal{Ann. Statist.}
\bvolume{31}
\bpages{768--806}.
\bid{doi={10.1214/aos/1056562462}, issn={0090-5364}, mr={1994730}}
\bptnote{check volume}%
\end{barticle}
%
\bptok{imsref}%
\endbibitem

\bibitem[\protect\citeauthoryear{Kline and Santos}{2012}]{Kline2012higher}
\begin{barticle}[mr]
\bauthor{\bsnm{Kline},~\bfnm{Patrick}\binits{P.}} \AND
\bauthor{\bsnm{Santos},~\bfnm{Andres}\binits{A.}}
(\byear{2012}).
\btitle{Higher order properties of the wild bootstrap under misspecification}.
\bjournal{J. Econometrics}
\bvolume{171}
\bpages{54--70}.
\bid{doi={10.1016/j.jeconom.2012.06.001}, issn={0304-4076}, mr={2970336}}
\end{barticle}
%
\bptok{imsref}%
\endbibitem

\bibitem[\protect\citeauthoryear{Koenker and Bassett}{1978}]{Koenker1978regression}
\begin{barticle}[mr]
\bauthor{\bsnm{Koenker},~\bfnm{Roger}\binits{R.}} \AND
\bauthor{\bsnm{Bassett},~\bfnm{Gilbert}\binits{G.} \bsuffix{Jr.}}
(\byear{1978}).
\btitle{Regression quantiles}.
\bjournal{Econometrica}
\bvolume{46}
\bpages{33--50}.
\bid{issn={0012-9682}, mr={0474644}}
\bptnote{check volume}%
\end{barticle}
%
\bptok{imsref}%
\endbibitem

\bibitem[\protect\citeauthoryear{Lavergne and Patilea}{2013}]{Lavergne2013smooth}
\begin{barticle}[mr]
\bauthor{\bsnm{Lavergne},~\bfnm{Pascal}\binits{P.}} \AND
\bauthor{\bsnm{Patilea},~\bfnm{Valentin}\binits{V.}}
(\byear{2013}).
\btitle{Smooth minimum distance estimation and testing with conditional estimating equations: Uniform in bandwidth theory}.
\bjournal{J. Econometrics}
\bvolume{177}
\bpages{47--59}.
\bid{doi={10.1016/j.jeconom.2013.05.006}, issn={0304-4076}, mr={3103915}}
\end{barticle}
%
\bptok{imsref}%
\endbibitem

\bibitem[\protect\citeauthoryear{Liu}{1988}]{Liu1988}
\begin{barticle}[mr]
\bauthor{\bsnm{Liu},~\bfnm{Regina~Y.}\binits{R.~Y.}}
(\byear{1988}).
\btitle{Bootstrap procedures under some non-i.i.d. models}.
\bjournal{Ann. Statist.}
\bvolume{16}
\bpages{1696--1708}.
\bid{doi={10.1214/aos/1176351062}, issn={0090-5364}, mr={0964947}}
\end{barticle}
%
\bptok{imsref}%
\endbibitem

\bibitem[\protect\citeauthoryear{Ma and Kosorok}{2005}]{MaKosorok2005robust}
\begin{barticle}[mr]
\bauthor{\bsnm{Ma},~\bfnm{Shuangge}\binits{S.}} \AND
\bauthor{\bsnm{Kosorok},~\bfnm{Michael~R.}\binits{M.~R.}}
(\byear{2005}).
\btitle{Robust semiparametric M-estimation and the weighted bootstrap}.
\bjournal{J. Multivariate Anal.}
\bvolume{96}
\bpages{190--217}.
\bid{doi={10.1016/j.jmva.2004.09.008}, issn={0047-259X}, mr={2202406}}
\end{barticle}
%
\bptok{imsref}%
\endbibitem

\bibitem[\protect\citeauthoryear{Mammen}{1992}]{Mammen1992does}
\begin{bbook}[auto:parserefs-M02]
\bauthor{\bsnm{Mammen},~\bfnm{E.}\binits{E.}}
(\byear{1992}).
\btitle{When Does Bootstrap Work?} \bseries{Lecture Notes in
Statistics}
\bvolume{77}.
\bpublisher{Springer},
\blocation{New York}.
\end{bbook}
%
\bptok{imsref}%
\endbibitem

\bibitem[\protect\citeauthoryear{Mammen}{1993}]{Mammen1993bootstrap}
\begin{barticle}[mr]
\bauthor{\bsnm{Mammen},~\bfnm{Enno}\binits{E.}}
(\byear{1993}).
\btitle{Bootstrap and wild bootstrap for high-dimensional linear models}.
\bjournal{Ann. Statist.}
\bvolume{21}
\bpages{255--285}.
\bid{doi={10.1214/aos/1176349025}, issn={0090-5364}, mr={1212176}}
\bptnote{check volume}%
\end{barticle}
%
\bptok{imsref}%
\endbibitem

\bibitem[\protect\citeauthoryear{Newton and Raftery}{1994}]{NewtonRaft1994WLB}
\begin{barticle}[mr]
\bauthor{\bsnm{Newton},~\bfnm{Michael~A.}\binits{M.~A.}} \AND
\bauthor{\bsnm{Raftery},~\bfnm{Adrian~E.}\binits{A.~E.}}
(\byear{1994}).
\btitle{Approximate {B}ayesian inference with the weighted likelihood bootstrap}.
\bjournal{J. R. Stat. Soc. Ser. B. Stat. Methodol.}
\bvolume{56}
\bpages{3--48}.
\bid{issn={0035-9246}, mr={1257793}}
\bptnote{check volume, check related}%
\end{barticle}
%
\bptok{imsref}%
\endbibitem

\bibitem[\protect\citeauthoryear{Spokoiny}{2012}]{Sp2012Pa}
\begin{barticle}[mr]
\bauthor{\bsnm{Spokoiny},~\bfnm{Vladimir}\binits{V.}}
(\byear{2012}).
\btitle{Parametric estimation. {F}inite sample theory}.
\bjournal{Ann. Statist.}
\bvolume{40}
\bpages{2877--2909}.
\bid{doi={10.1214/12-AOS1054}, issn={0090-5364}, mr={3097963}}
\end{barticle}
%
\bptok{imsref}%
\endbibitem

\bibitem[\protect\citeauthoryear{Spokoiny}{2013}]{Spokoiny2013Bernstein}
\begin{bmisc}[auto:parserefs-M02]
\bauthor{\bsnm{Spokoiny},~\bfnm{V.}\binits{V.}}
(\byear{2013}).
\bhowpublished{Bernstein--von {M}ises {t}heorem for growing parameter dimension.
Preprint. Available at \arxivurl{arXiv:1302.3430}}.
\end{bmisc}
%
\bptok{imsref}%
\endbibitem


\bibitem[\protect\citeauthoryear{Spokoiny and Zhilova}{2015}]{SpZh2014PMBsupp}
\begin{bmisc}[author]
\bauthor{\bsnm{Spokoiny},~\bfnm{V.}\binits{V.}} \AND
\bauthor{\bsnm{Zhilova},~\bfnm{M.}\binits{M.}}
(\byear{2015}).
\bhowpublished{Supplement to ``Bootstrap confidence sets under model misspecification.''
DOI:\doiurl{10.1214/15-AOS1355SUPP}}.
\bptok{imsref}%
\end{bmisc}
\endbibitem
\bptok{imsref}%
\endbibitem

\bibitem[\protect\citeauthoryear{van~der Vaart and Wellner}{1996}]{VaartWellner1996weak}
\begin{bbook}[mr]
\bauthor{\bsnm{van~der Vaart},~\bfnm{Aad~W.}\binits{A.~W.}} \AND
\bauthor{\bsnm{Wellner},~\bfnm{Jon~A.}\binits{J.~A.}}
(\byear{1996}).
\btitle{Weak Convergence and Empirical Processes}.
\bpublisher{Springer},
\blocation{New York}.
\bid{doi={10.1007/978-1-4757-2545-2}, mr={1385671}}
\end{bbook}
%
\bptok{imsref}%
\endbibitem

\bibitem[\protect\citeauthoryear{Wilks}{1938}]{wilks1938}
\begin{barticle}[auto:parserefs-M02]
\bauthor{\bsnm{Wilks},~\bfnm{S.~S.}\binits{S.~S.}}
(\byear{1938}).
\btitle{The large-sample distribution of the likelihood ratio for testing composite hypotheses}.
\bjournal{The Annals of Mathematical Statistics}
\bvolume{9}
\bpages{60--62}.
\end{barticle}
%
\bptok{imsref}%
\endbibitem

\bibitem[\protect\citeauthoryear{Wu}{1986}]{Wu1986wildboot}
\begin{barticle}[mr]
\bauthor{\bsnm{Wu},~\bfnm{C.-F.~J.}\binits{C.-F.~J.}}
(\byear{1986}).
\btitle{Jackknife, bootstrap and other resampling methods in regression analysis}.
\bjournal{Ann. Statist.}
\bvolume{14}
\bpages{1261--1350}.
\bid{doi={10.1214/aos/1176350142}, issn={0090-5364}, mr={0868303}}
\bptnote{check related, check pages}%
\end{barticle}
%
\bptok{imsref}%
\endbibitem
\end{thebibliography}




\printaddresses
\end{document}